\documentclass[UTF8]{article}

\usepackage{subcaption}
\usepackage{algorithm}
\usepackage{algorithmicx}
\usepackage{algpseudocode}
\usepackage{amsmath}
\usepackage{amsfonts}
\usepackage{amssymb}
\usepackage{bm}
\usepackage{color}
\usepackage{array,multirow,booktabs,bbm} 
\usepackage{enumitem}
\usepackage{epstopdf}
\usepackage{float}
\usepackage{graphicx}
\usepackage{indentfirst}
\usepackage{listings}
\newcommand{\Rmnum}[1]{\expandafter\@slowromancap\romannumeral #1@}
\usepackage{makeidx}
\usepackage{multicol}
\usepackage{mathrsfs}
\usepackage{mathtools}
\usepackage{soul}
\usepackage{tikz}
\usepackage{verbatim}
\usepackage{xcolor}
\usepackage{verbatim}
\usepackage{amsthm}
\usepackage{tabularx}
\usepackage{threeparttable,multirow,dcolumn,booktabs}

\newcommand{\bs}[1]{\boldsymbol{#1}}
\newcommand{\bmu}{\bs{\mu}}

\newcommand{\bff}{\bs{f}}
\newcommand{\bu}{\bs{u}}
\newcommand{\br}{\bs{r}}
\newcommand{\bss}{\bs{s}}

\usepackage[normalem]{ulem}

\newcommand{\bx}{\bs{x}}

\newcommand{\bV}{\bs{V}}
\newcommand{\bH}{\bs{H}}

\newcommand{\bE}{\bs{E}}
\newcommand{\bB}{\bs{B}}

\newcommand{\bxi}{\bs{\xi}}

\newcommand{\bS}{\bs{S}}
\newcommand{\bA}{\bs{A}}
\newcommand{\tr}{\tilde{{r}}}
\newcommand{\bW}{\bs{W}}

\newcommand{\s}{\bs{s}_r}

\newcommand{\calN}{\mathcal{N}}

\DeclareMathOperator*{\argmax}{argmax}
\DeclareMathOperator*{\argmin}{argmin}

\title{Nonlinear Model Order Reduction on Quadratic Manifolds via Greedy Algorithms with Dimension-Dependent Regularization}

\author{
{Lijie Ji}\thanks{Department of Mathematics, Shanghai University, Shanghai 200444, P.R. China; Newtouch Center for Mathematics of Shanghai University, Shanghai University, Shanghai 200444, P.R. China. Email: {\tt lijieji@shu.edu.cn}.},
\quad Sabrina Rashid \footnote{Department of Mathematics, University of South Carolina, Columbia, SC 29208, USA. Email: {\tt{SR89@email.sc.edu}}.}, \quad
Yanlai Chen\footnote{Department of Mathematics, University of Massachusetts Dartmouth, 285 Old Westport Road, North Dartmouth, MA 02747, USA. Email: {\tt{yanlai.chen@umassd.edu}}. Y. Chen was partially supported by National Science Foundation grant DMS-2208277},\quad
Zhu Wang\footnote{
Department of Mathematics, University of South Carolina, Columbia, SC 29208, USA. Email: {\tt{wangzhu@math.sc.edu}}. Z. Wang was partially supported by the U.S. National Science Foundation under Award Nos. DMS-2038080, DMS-2245097, and DMS-2513112.}
}

\makeatother

\begin{document}
\maketitle

\begin{abstract}
Traditional projection-based reduced-order modeling approximates the full-order model by projecting it onto a linear subspace. With a fast-decaying Kolmogorov $n$-width of the solution manifold, the resulting reduced-order model (ROM) can be an efficient and accurate emulator. However, for parametric partial differential equations with slowly decaying Kolmogorov $n$-width, the dimension of the linear subspace required for a reasonable accuracy becomes very large, which undermines computational efficiency. To address this limitation, quadratic manifold methods have recently been proposed. These data-driven methods first identify a quadratic mapping by minimizing the linear projection error over a large set of snapshots, often with the aid of regularization techniques to solve the associated minimization problem, and then use this mapping to construct ROMs.

In this paper, we propose and test a novel enhancement to this quadratic manifold approach by introducing a first-of-its-kind double-greedy algorithm on the regularization parameters coupled with a standard greedy algorithm on the physical parameter. Our approach balances the trade-off between the accuracy of the quadratic mapping and the stability of the resulting nonlinear ROM, leading to a highly efficient and data-sparse algorithm. 
Numerical experiments conducted on equations such as linear transport, acoustic wave,  advection-diffusion, and Burgers' demonstrate the accuracy, efficiency, and stability of the proposed algorithm.

\end{abstract}
\section{Introduction}

Reduced-order models (ROMs) have been developed for various multi-query and real-time applications, including uncertainty quantification \cite{resseguier2021quantifying,chen2017reduced}, parametric computational fluid dynamics \cite{stabile2018finite,Kunisch_Volkwein_POD,wang2012proper}, and optimal control \cite{negri2013reduced,karcher2018certified,bader2016certified}. The general idea is to first identify a low-dimensional approximation of the solution manifold associated with the governing equations, and then to construct a low-dimensional surrogate of the original full-order model (FOM). One popular class of reduced-order modeling approaches is linear projection-based ROMs \cite{benner2015survey,machiels2001output,chen2010certified}, in which {case a linear basis matrix $\bV_r \in \mathbb{R}^{\calN \times r}$ is constructed from solution data, providing a low-dimensional linear subspace spanned by the columns of $\bV_r$}. The state variable $\bu \in \mathbb{R}^\calN$ is then approximated by a linear combination of these basis vectors, that is, $\bu \approx \bV_r \s$, where $\s \in \mathbb{R}^{r}$ denotes the reduced state. Such {linear} ROMs, {with ``linear" referring} to the reconstruction map rather than the governing equations, can be built either intrusively, via Galerkin or Petrov-Galerkin projection of the governing equations \cite{Quarteroni2015,HesthavenRozzaStammBook,CarlbergBouMoslehFarhat2011}, or non-intrusively, through interpolation or regression methods \cite{kramer2024learning,xiao2016non,HesthavenUbbiali,guo2018reduced,CHEN2021110666}. However, for certain types of problems such as transport-dominated problems, it is challenging to obtain accurate approximations  using linear ROMs. From a mathematical perspective, the Kolmogorov $n$-width of the associated solution manifold decays slowly in these cases, indicating that a large number of basis vectors is required to achieve a prescribed accuracy \cite{cohen2016kolmogorov}. As a result, linear ROMs often struggle to provide efficient and accurate reduced-order approximations for such problems \cite{greif2019decay,cohen2023nonlinear}.

{Given the limitations of linear reconstruction map, nonlinear model reduction techniques have been developed with the goal of constructing a low-dimensional surrogate via a nonlinear reconstruction map.} In this work, we focus on the nonlinear ROM built on a quadratic manifold (QM) \cite{barnett2022quadratic}, where the approximated solution takes the form $\widehat{\bu} = \bV_r \s + \bH_r (\s \widetilde{\otimes} \s)$. Here, $\bH_r \in \mathbb{R}^{\calN \times q}$ {represents the quadratic mapping operator in matrix form with $q = r(r+1)/2$}, and $\widetilde{\otimes}$ denotes the modified Kronecker product that extracts only the unique quadratic terms from the standard Kronecker product $\s \otimes \s$, yielding a compressed vector of dimension $q$. The conventional QM approach typically constructs $\bV_r$ and $\bH_r$ in a purely data-driven manner, which requires {many queries} of the FOM. Specifically, high-fidelity solutions corresponding to the training parameters are collected, and snapshots (sampled at every \( l_{\text{sam}}\)-th time step) are assembled into a snapshot matrix. A singular value decomposition (SVD) of this matrix is performed, and $\bV_r$ is formed by selecting the first \( r \) left singular vectors. Subsequently, $\bH_r$ is derived by minimizing the linear projection error of the snapshot matrix onto $\bV_r$. Using the nonlinear reconstruction form in QM, one can derive the ROM via Galerkin projection. In particular, the stability and accuracy of the resulting ROM depend on the regularization parameter used in the optimization method, such as the Tikhonov regularization method. In  \cite{barnett2022quadratic,schwerdtner2024greedy}, the regularization parameter is determined by minimizing the reconstruction error on a validation set or through generalized cross-validation. More details on the construction of $\bH_r$ are provided in Section~\ref{sec:qm}.

{In this work, we propose a novel ROM constructed via greedy algorithms with dimension-dependent regularization on a quadratic manifold} (Greedy-Quadratic ROM). The main contributions are as follows:
\begin{itemize}
    \item \textbf{Greedy construction of the reduced basis matrix:} Instead of generating both linear basis matrix and quadratic mapping matrix from a large set of snapshots, a residual-based \textit{a posteriori} error estimator is designed to incrementally select representative parameters and hierarchically construct the linear basis matrix, after which the quadratic mapping matrix is learned in a data-driven manner.
    \item \textbf{A stabilized quadratic manifold:} To balance the stability of the ROM and the accuracy of the QM, a double-greedy algorithm is developed to choose the regularization parameters used in  generating the quadratic mapping matrix. Consequently, it establishes a regularization mechanism that depends on the dimension of the linear basis, $r$, and ensures alignment between the reconstruction error and the ROM error.

\end{itemize}
Although greedy frameworks have been widely used in linear ROMs \cite{HesthavenRozzaStammBook, Quarteroni2015}, to the best of our knowledge, this is the first work that integrates a greedy framework with the construction of a quadratic mapping matrix while simultaneously addressing both stability and accuracy of the reduced model. By incorporating hyper-reduction techniques such as the empirical interpolation method \cite{Barrault2004,klein2023energy}, both the offline and online computational efficiency of the proposed algorithm can be further enhanced.

\subsection{{A motivating example}}

{We demonstrate the need of our double-greedy algorithm by showcasing the delicate balance between stability and accuracy, and their sensitivity to the regularization parameter via} the one-dimensional linear transport equation {parameterized by location of the Gaussian pulse in the initial condition}:
\begin{align}
\frac{\partial}{\partial t} u(x, t) & +c \frac{\partial}{\partial x} u(x, t)=0, \quad x \in [0,1], t \in(0, T),\\
u_0(x) & :=u(x, 0)=\frac{1}{\sqrt{2\pi}\,\sigma} \exp\left(-\frac{(x-x^\ast)^2}{2\sigma^2}\right). \nonumber
\label{eq:intro:transport}
\end{align}
Here, \( c \) denotes the constant advection velocity. 
The high-fidelity numerical solution is obtained by discretizing the time derivative \( \partial u / \partial t \) with the Lax--Wendroff scheme and the spatial derivative with a central finite-difference scheme. The spatial domain is uniformly discretized into \( \calN = 2000 \) grid points, and the time domain is divided into \( {N_T} = 4000 \) uniform steps. The parameters are set as \( c = 10 \), \( \sigma = 0.01 \), and the peak position \( x^\ast \) of the initial Gaussian pulse is varied within the interval \( [0.05,\,0.25] \). The training set is generated by uniformly sampling \( 41 \) points from \( [0.05,0.25] \), while the testing set consists of \( 5 \) uniformly distributed points taken from the sub-interval \( [0.0524,\,0.226] \).

\begin{figure}[htb!]
\centering
\includegraphics[scale=0.24]{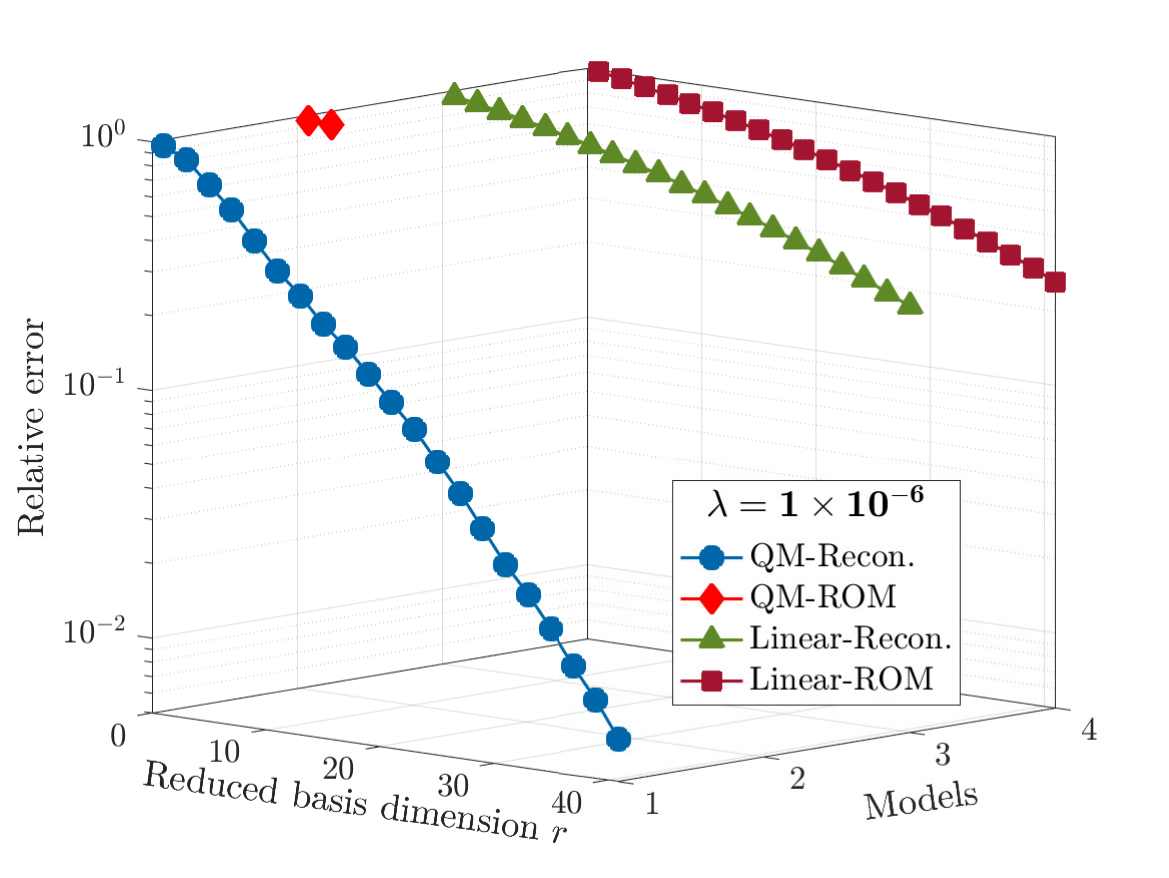}
\includegraphics[scale=0.24]{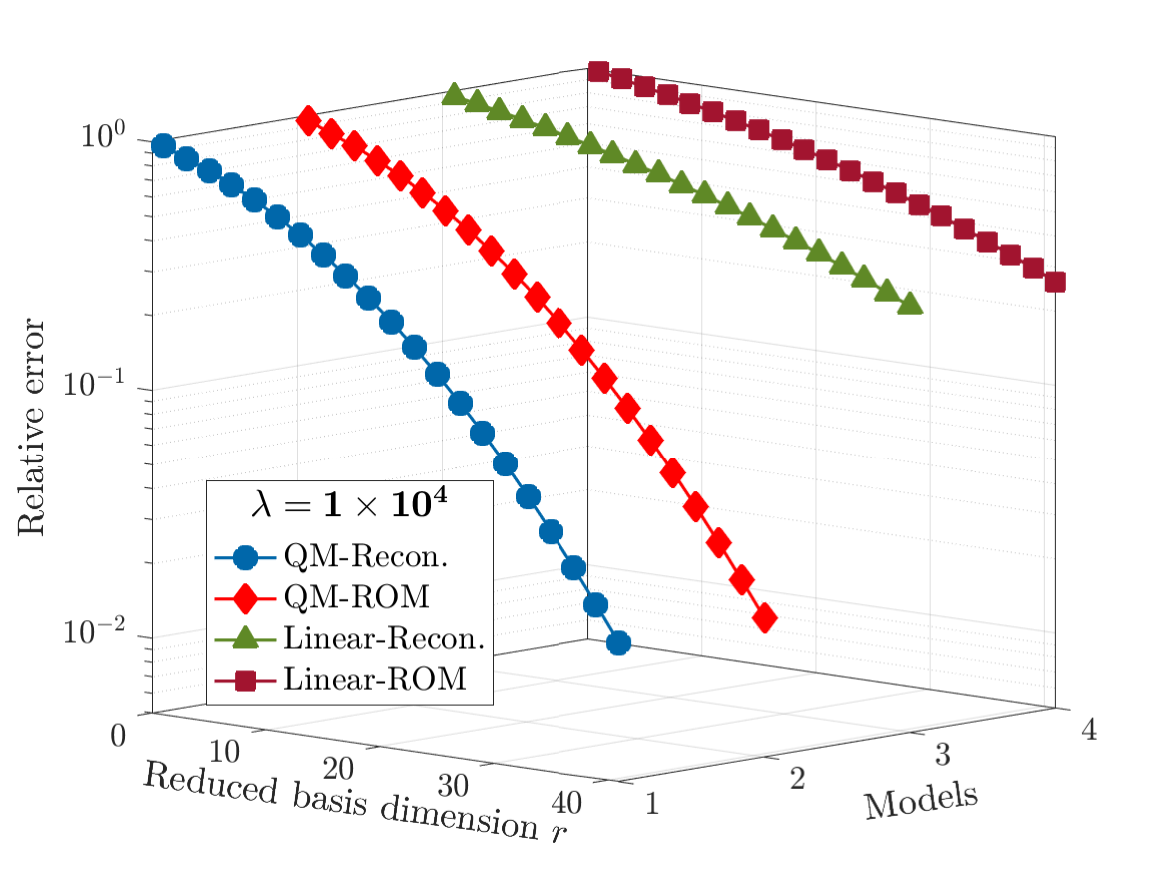}
\includegraphics[scale=0.24]{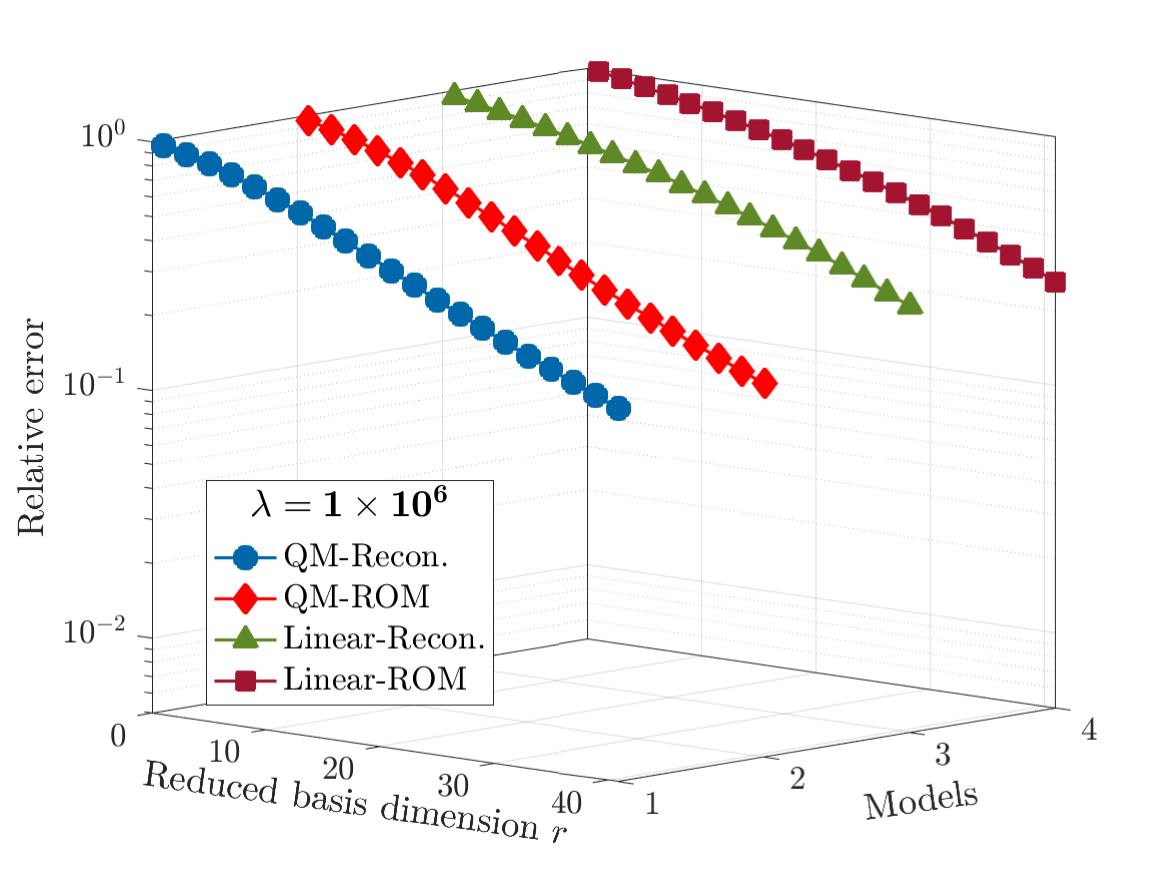}
\caption{The relative linear-/QM-ROM and reconstruction errors as a function of reduced basis dimensions $r$, for three regularization parameters $\lambda=1e-6, 1e4, 1e6$ (from Left to Right). For $\lambda=1e-6$, the QM-ROM is unstable for $r\geq 5$, producing infinite relative errors that are not shown in the left panel. For $\lambda=1e6$ in the right panel, the QM-ROM barely outperforms the linear-ROM. Moderate $\lambda=1e4$ provides the best balance between stability and accuracy among these three cases.}
\label{fig:intro:1}
\end{figure}

To illustrate the influence of the regularization parameter, Figure~\ref{fig:intro:1} displays the average relative errors of the linear- and QM-ROMs (and their corresponding reconstructions {errors}) on the testing set for different choices of the regularization parameter. Here, ``Linear-/QM-ROM'' refers to the error obtained when the reduced state is computed by solving the ROM under the Galerkin projection, whereas ``Linear-/QM-Recon.'' denotes the error resulting from a direct projection \( \s = \mathbf{V}_r^{\mathsf{T}} \bu \). The precise definition of the relative error is given in Section~\ref{sec:numerical}. As shown in the left panel, when the regularization parameter $\lambda=10^{-6}$ is small, the quadratic approximation provides a significantly more accurate reconstruction than the linear approximation for large values of $r$. {However, in this case, the QM-ROM is stable only for very low reduced dimensions, namely $r = 1$ and $3$. For $r \geq 5$, the QM-ROM becomes unstable (with corresponding infinite relative errors not shown in the plot).} When $\lambda$ is increased to $10^4$, as shown in the middle panel, the QM-ROM approximation accuracy coincides with that of the QM reconstruction; however, the reconstruction accuracy is worse than that obtained for $\lambda =10^{-6}$. When $\lambda=10^6$, as shown in the right panel, the QM-ROM approximation accuracy also coincides with that of the QM reconstruction; however, this accuracy is only slightly better than the linear ones. These observations indicate that (i) the choice of the regularization parameter, $\lambda$, is critical for the numerical stability of the constructed QM-ROMs, and (ii) stability and accuracy must be balanced by selecting an appropriate amount of regularization.

\subsection{Literature review}
To overcome the central challenge of a slowly decaying Kolmogorov $n$-width (a phenomenon that renders traditional ROM based on a single global linear subspace inefficient for a broad class of problems, from nonlinear multiphysics systems to linear hyperbolic transport equations), researchers have developed diverse modeling strategies. Among these, one prominent category employs localization, such as by partitioning the time, parameter, or physical domains \cite{amsallem2012nonlinear,San2015PrincipalID,San2018}.
For instance, the principal interval decomposition (PID) is applied to the time domain for incompressible Navier–Stokes equations in \cite{klein2023energy}, where an energy-conserving hyper-reduction technique further accelerates online computation. In \cite{cheung2023local}, a Lagrangian-based ROM with temporal domain partitioning is introduced for the Rayleigh–Taylor instability problem. Meanwhile, researchers in \cite{haasdonk2011training} propose a grid-based adaptive partitioning of the parameter domain, constructing multiple local bases tailored to different parameter regions.
Spatial localization has also been extensively explored. A localized discrete empirical interpolation method (L-DEIM) is introduced in \cite{PeherstorferButnaruWillcoxBungartz2014}, where the nonlinear term is approximated using multiple local basis patches rather than a single global basis. Further studies on spatially local ROMs via physical domain decomposition can be found in \cite{PAN2024113452, XIAO201915}, among other related works.

The second category involves the use of an online adaptive basis matrix rather than a fixed basis matrix \cite{peherstorfer2015online}. The sampling points in the associated hyper-reduction method are further enhanced using residual information in \cite{Benjamin2020}. The third prominent model reduction techniques include learning nonlinear encoder-decoder mappings or employing parametrized coordinate transformations. For instance, \cite{reiss2018shifted} proposed a shifted proper orthogonal decomposition (sPOD), which remains effective for problems involving multiple transport velocities. A data-driven neural network variant of sPOD was developed in \cite{NNsPOD2022} and validated in multiphase flow simulations. Registration methods learn domain mappings which tracks moving features and a linear reduction is performed in transformed snapshots \cite{mirhoseini2023model, barral2023registration, ferrero2022registration, taddei2020registration, PENG2023114560}.  Encoder-decoder frameworks based on deep learning have also attracted significant research interest \cite{duan2024non, lee2020model}. Alternatively, some approaches reformulate the original parametric PDE in its Lagrangian form and subsequently apply linear reduction techniques, as explored in \cite{mojgani2017lagrangian, LU2020109229}.

In addition to the methods mentioned above, the QM approach has gained attention in recent years. This method represents the state as a linear reduced approximation augmented by a quadratic correction term \cite{yong2025learning}, an idea originally proposed for structural dynamics \cite{jain2017quadratic} and later generalized within a data-driven framework \cite{barnett2022quadratic}. Subsequent developments have integrated it with operator inference for non-intrusive computation \cite{geelen2023operator} and augmented it with neural networks for greater flexibility \cite{barnett2023neural}.
A central challenge in this approach is the construction of an optimal linear subspace. The straightforward use of the leading POD modes often proves suboptimal. To address this, recent research has focused on co-optimizing the linear and nonlinear components. For instance, \cite{schwerdtner2024greedy} proposed a greedy algorithm that iteratively selects singular vectors to minimize the reconstruction error of the QM. Similarly, \cite{geelen2024learning} introduced an alternating optimization scheme to jointly learn the linear and the quadratic mapping.
Alongside the linear and nonlinear mapping decomposition framework, the authors of \cite{guo2025probabilistic} employ POD for linear compression but construct the nonlinear approximation via a probabilistic manifold, demonstrating its efficacy in fluid dynamics problems.

\subsection{Organizations}
The remainder of this paper is organized as follows. In Section \ref{sec:qm}, the general parametrized time-dependent PDE is introduced, and the nonlinear ROM with quadratic manifold approximation is reviewed. Section \ref{sec:greedyqm} presents the ROM with greedy algorithms on a quadratic manifold, including the a \textit{posteriori} error estimator, and the double-greedy algorithm. Extensive numerical tests in Section \ref{sec:numerical} demonstrate the accuracy, stability, and efficiency of the proposed method. Conclusions are provided in Section \ref{sec:conclusion}.

\section{Nonlinear ROM on a quadratic manifold}
\label{sec:qm}

\subsection{Notation and model equation}
\label{sec:notation}
For clarity, we summarize the primary mathematical notation used throughout this paper in Table~\ref{tab:notation}. The conventions for vectors and matrices are as follows: bold lowercase letters (e.g., $\bu$, $\s$) denote vectors, bold uppercase letters (e.g., $\bV_r$, $\bA_i$) denote matrices.

For simplicity, we consider a linear time-dependent PDE with {affine} parametric dependence 
\begin{align}
\frac{\partial u}{\partial t} = \sum_{i=1}^{Q_a} \alpha_i(\bmu) L_i(u) + f(\bx), \quad \bx \in \Omega_{\bx},\ t \in [0, T],\ \bmu \in \Omega_p.
\end{align}
Here, $L_i$ denotes the $i$-th linear differential operator, $\alpha_i(\bmu)$ is the parameter-dependent coefficient, and $f$ is a parameter-independent function.
Using a proper spatial discretization, we derive the semi-discrete scheme as follows:
\begin{align}
\frac{\partial \bu}{\partial t} = \sum_{i=1}^{Q_a} \alpha_i(\bmu) \bA_i \bu + \bff,
\label{eq:fom}
\end{align}
where $\bu(t, \bmu) \in \mathbb{R}^{\calN}$ represents the numerical solution at time $t$, and $\calN$ is the total number of degrees of freedom in physical space. The matrix $\bA_i$ {is the finite discretization of} the linear operator $L_i$.

\begin{table}[htbp]
    \centering
    \caption{Some mathematical notations adopted in the paper.}
    \label{tab:notation}
    \renewcommand{\arraystretch}{1.2}
    \begin{tabularx}{\linewidth}{p{1.8cm} X}
        \toprule
        \textbf{Symbol} & \textbf{Description} \\
        \midrule
        $\bx$ & Spatial coordinate, belonging to the physical domain $\Omega_{\bx}$. \\
        $\bmu$ & Parameter vector, belonging to the parameter domain $\Omega_p$. \\
        $\bu$ & FOM state vector (spatially discretized solution) \\
        $\calN$ & Number of degrees of freedom (FOM dimension). \\
        $\bV_r$ & Linear basis matrix (orthonormal columns). \\
        $r$ & Dimension of the reduced linear subspace (ROM dimension). \\
        $\s$ & Reduced state vector. \\
        $\bH_r$ & Quadratic mapping matrix for nonlinear reconstruction. \\
        $q$ & Number of unique quadratic terms. \\
        $\widehat{\bu}$ & ROM approximation of the FOM state via the quadratic manifold. \\
        $N_T$ & Number of time steps in the FOM simulation. \\
        $\Xi_{\mathrm{train}}$ & Discrete training parameter set (a subset of $\Omega_p$). \\
        $n_{\text{train}}$ & Number of parameters in the training set $\Xi_{\mathrm{train}}$. \\
        $\bS$ & Snapshot matrix collecting sampled FOM solutions. \\
        $l$ & Total number of snapshots in $\bS$. \\
        $\bE$ & Linear projection (POD) error matrix. \\
        $\mathbf{W}$ & Matrix of reduced-state quadratic terms. \\
        $\lambda$ & Tikhonov regularization parameter. \\
        $E_{\lambda}$ & Frobenius norm residual associated with regularization parameter $\lambda$. \\
        $\Xi_{\lambda}$ & Discrete candidate set for the regularization parameter $\lambda$. \\
        $l_{\mathrm{sam}}$ & Temporal sampling stride (store snapshot every $l_{\mathrm{sam}}$ time steps). \\
        $\Delta_{r,m}(\bmu)$ & Residual-based \textit{a posteriori} error estimator for a given $\bmu$. \\
        $N_{\text{incre}}$ & Basis increment size in the POD-Greedy algorithm. \\
        \bottomrule
    \end{tabularx}
\end{table}

\subsection{Quadratic manifold as a linear-nonlinear autoencoder}

As stated in the introduction, in the quadratic manifold approach \cite{barnett2022quadratic}, the reduced state of the FOM at time $t$ is computed using a linear encoder once the linear basis matrix $\bV_r \in \mathbb{R}^{\calN \times r}$ is generated:
\begin{equation}
\s(t; \bmu) = \bV_r^\top \bu(t; \bmu),
\end{equation}
where $r$ is the dimension of the (linear) ROM. In contrast, the FOM solution is then approximated via a nonlinear decoder employing a quadratic manifold:
\begin{equation}
\bu(t; \bmu) \approx \widehat{\bu}(t; \bmu) = \bV_r \s(t;\bmu) + \bH_r (\s (t;\bmu)\widetilde{\otimes} \s (t;\bmu)).
\label{eq:decoder}
\end{equation}

\subsubsection*{Regularized generation of the quadratic mapping matrix}
The quadratic mapping matrix $\bH_r$ is typically determined in a purely data-driven manner. This procedure is detailed as follows. First, using a suitably designed time discretization scheme, high-fidelity solutions are obtained at all time instances: $\bu^j(\bmu),\ j = 0, \cdots, N_T$, where $N_T = T / \Delta t$ denotes the number of time steps, and $\Delta t$ is the time step size. After computing all high-fidelity solutions for parameters in the training set, they are arranged into a snapshot matrix:
\[
\bS = \left[ \bu^0(\bmu^1),\ \bu^{l_{\text{sam}}}(\bmu^1),\ \bu^{2l_{\text{sam}}}(\bmu^1)\, \cdots,\ \bu^{N_T}(\bmu^1),\ \cdots,\ \bu^{N_T}(\bmu^{n_{\text{train}}}) \right] \in \mathbb{R}^{\calN \times l},
\]
where $l = n_{\text{train}} \cdot (N_T / l_{\text{sam}}+1)$. Here, $n_{\text{train}}$ denotes the number of parameters in the training set and $l_{\text{sam}}$ is the sampling time step size.
In linear ROMs, the snapshot matrix $\bS$ is decomposed using SVD, and the first $r$ left singular vectors are selected to form the linear basis matrix $\bV_r$. The corresponding POD projection error (linear recovery error) is then given by:
\begin{align}
\bE = \bS - \bV_r \bV_r^\top \bS = (I - \bV_r \bV_r^\top) \bS.
\label{eq:poderror}
\end{align}
The quadratic mapping matrix is designed to minimize this error by leveraging quadratic information of the reduced state $\s$:
\begin{align}
\bH_r = \underset{\bH_r \in \mathbb{R}^{\calN \times q}}{\arg \min } \frac{1}{2} \left\| \mathbf{W}^{\top} \bH_r^{\top} - \bE^{\top} \right\|_F^2,
\label{eq:qm_H_noreg}
\end{align}
where
\begin{align}
\mathbf{W} := \left( \begin{array}{cccc}
\mid & \mid & & \mid \\
\bss_r^1 \widetilde{\otimes} \bss_r^1 & \bss_r^2 \widetilde{\otimes} \bss_r^2 & \ldots & \bss_r^l \widetilde{\otimes} \bss_r^l \\
\mid & \mid & & \mid
\end{array} \right) \in \mathbb{R}^{q \times l}.
\label{eq:WMatrix}
\end{align}
This optimization problem often requires regularization. In this work, we employ Tikhonov regularization \cite{hansen1998rank} {converting \eqref{eq:qm_H_noreg} to its regularized version}:
\begin{align}
\bH_r = \underset{\bH_r \in \mathbb{R}^{\calN \times q}}{\arg \min } \frac{1}{2} \left\| \mathbf{W}^{\top} \bH_r^{\top} - \bE^{\top} \right\|_F^2 + \frac{\lambda}{2} \left\| \bH_r^{\top} \right\|_F^2.
\label{eq:qm:H}
\end{align} 
In practical calculations, each row of $\bH_r$ can be solved independently. Denoting by $\mathbf{h}_i$ the $i$-th row of $\bH_r$, we solve $\calN$ independent overdetermined problems:
\[
\mathbf{h}_i = \arg\min _{\beta \in \mathbb{R}^{ q \times 1}} \frac{1}{2}\left\| \bE^{\top}(:, i) - \mathbf{W}^{\top} \beta \right\|_2^2 + \frac{\lambda}{2} \left\| \beta \right\|_2^2, \quad i = 1, \ldots, \calN,
\]
where $\mathbf{h}_i = \bH_r^{\top}(:, i)$. These least-squares problems admit explicit solutions, and all the $\mathbf{h}_i$ can be calculated simultaneously; see Algorithm \ref{alg:tikhonov}. 
\begin{algorithm}[h]
\caption{Tikhonov Regularized Least-Squares Solution for the Quadratic Mapping Matrix}
\label{alg:tikhonov}
\begin{algorithmic}[1]
    \Require Matrix $\mathbf{A}:=\mathbf{W}^{\top}$, matrix $\mathbf{B}:= \bE^{\top}$, regularization parameter $\lambda$
    \Ensure Quadratic mapping matrix $\bH_r$, residual norm $E_{\lambda}$
    \Statex
    \State $[\mathbf{U}, \Sigma, \mathbf{V}] \gets \operatorname{svd}(\mathbf{A})$ \Comment{Compute SVD}
    \State $\boldsymbol{\sigma} \gets \operatorname{diag}(\Sigma)$ \Comment{Extract singular values}
    \State $\mathbf{s}_{\text{inv}} \gets \boldsymbol{\sigma} / (\boldsymbol{\sigma}^2 + \lambda^2)$ \Comment{Compute Tikhonov filter factors}
    \State $\mathbf{X} \gets \mathbf{V} \cdot \operatorname{diag}(\mathbf{s}_{\text{inv}}) \cdot (\mathbf{U}^\top \mathbf{B})$ \Comment{Compute regularized solution}
    \State $E_\lambda \gets \|\mathbf{B} - \mathbf{A} \mathbf{X}\|_F$ \Comment{Compute residual norm}
    \State $\bH_r \gets \mathbf{X}^\top$ \Comment{Transpose to obtain $\bH_r$}
    \State \Return $\bH_r, E_{\lambda}$
\end{algorithmic}
\end{algorithm}

\subsubsection*{Reduced order formulation for quadratic manifold}
The quadratic mapping matrix $\bH_r$ is constructed by minimizing the projection error of the linear basis matrix $\bV_r$; consequently, $\bH_r$ is orthogonal to $\bV_r$.
Substituting the nonlinear decoder expression from Eq. \eqref{eq:decoder} into the semi-discrete formulation yields
\begin{align}
\frac{\partial [ \bV_r \s(t) + \bH_r(\s \widetilde{\otimes} \s)]}{\partial t} = \sum_{i=1}^{Q_a} \alpha_i(\bmu) \bA_i [\bV_r \s(t) + \bH_r(\s \widetilde{\otimes} \s)] + \bff.
\end{align}
Applying Galerkin projection and utilizing the orthogonality between $\bV_r$ and $\bH_r$ {and the fact that columns of $\bV_r$ are orthonormal, we obtain} the following semi-discrete reduced system
\begin{align}
\frac{\partial \s(t)}{\partial t} = \sum_{i=1}^{Q_a} \alpha_i(\bmu) \left[ \widehat{\bA}_i \s(t) + \widehat{\bB}_i (\s \widetilde{\otimes} \s) \right] + \widehat{\bff},
\label{eq:galerkinrom}
\end{align}
where $\widehat{\bA}_i = \bV_r^\top \bA_i \bV_r \in \mathbb{R}^{r \times r}$, $\widehat{\bB}_{i} = \bV_r^\top \bA_i \bH_r \in \mathbb{R}^{r \times \tr}$, and $\widehat{\bff} = \bV_r^\top \bff \in \mathbb{R}^{r \times 1}$. 
Since $r \ll \calN$ typically holds, solving this reduced system is computationally significantly cheaper than solving the FOM in Eq. \eqref{eq:fom}. However, it should be noted that the incorporation of quadratic terms in the reduced state introduces additional nonlinearity to the reduced model, which consequently requires careful consideration of numerical stability properties, as indicated in Figure \ref{fig:intro:1}. Thus, choosing an appropriate quadratic mapping matrix, that is, selecting a suitable regularization parameter, is curial for balancing the accuracy and stability of the ROM in Eq. \eqref{eq:galerkinrom}.

\section{Greedy algorithms on quadratic manifold}
\label{sec:greedyqm}

The quadratic mapping matrix is conventionally obtained through data-driven methods that rely on high-fidelity solutions computed for extensive parameter samplings, as shown in Section \ref{sec:qm}. This computational approach requires numerous accesses to the FOM, resulting in high computational costs. As known, an \textit{a posteriori} error estimator is often designed in reduced basis method, and only high-fidelity solutions for selected representative parameters are calculated in the offline process \cite{Quarteroni2015}.
In this section, we introduce our ROM on a quadratic manifold, which includes three key components: an \textit{a posteriori} error estimator for representative parameters' selection, the POD-Greedy algorithm for generating the linear basis matrix, and a double-greedy algorithm for determining the regularization parameter during the generation of the quadratic mapping matrix.

\textbf{\textit{A posteriori} error estimator:}
Given the reduced states \(\{\s^{j}\}_{j=0}^{N_T}\) obtained by discretizing (e.g., via the explicit Euler method) and solving the reduced system \eqref{eq:galerkinrom}, we reconstruct the corresponding full-order approximations at times \(t_j\) and \(t_{j-1}\) using the quadratic manifold decoder \eqref{eq:decoder}:
\[
\widehat{\bu}^{j} = \bV_r \s^{j} + \bH_r (\s^{j} \widetilde{\otimes} \s^{j}), \quad j=0, \dots, N_T.
\]
We then define the residual vector \(\br^{j} \in \mathbb{R}^{\calN}\) at time \(t_j\) as the amount by which the reconstructed solution \(\widehat{\bu}^{\,j}\) fails to satisfy the discretized full-order equations:
\begin{align}
\br^j = \widehat{\bu}^j -\widehat{\bu}^{j-1} -\Delta t\sum_{i=1}^{Q_a} \alpha_i(\bmu) \bA_i \widehat{\bu}^{j-1} -\Delta t\bff, \quad j=1, 2, \cdots, N_T.
\end{align}
We finally formulate the residual-based \textit{a posteriori} error estimator for the ROM with an \(r\)-dimensional basis matrix at the \(m\)-th greedy iteration as:
\begin{align}
\Delta_{r,m} (\bmu) = \sqrt{\sum_{j=1}^{N_T} \lVert \br^j \rVert_2^2}.
\end{align}
The $(m+1)$-th parameter is then determined by maximizing this \textit{a posteriori} error estimator over a discrete training set $\Xi_{\rm train} \subset \Omega_p$:
\begin{align}
\bmu^{m+1} = \argmax_{\bmu \in \Xi_{\rm train}} \Delta_{r,m}(\bmu).
\label{eq:errorestimator}
\end{align}

\textbf{POD-Greedy algorithm:} 
In the linear ROM framework, the POD-Greedy method is widely used for constructing the linear basis matrix in time-dependent parametrized partial differential equations \cite{hesthaven2022reduced,haasdonk2008reduced,haasdonk2013convergence}. After the $(m+1)$-th representative parameter is determined by maximizing the residual-based \textit{a posteriori} error estimator, its high-fidelity solution is computed. The projection error onto the current linear reduced basis matrix $\bV_r$ is then calculated. The basis matrix is subsequently enriched by extracting the first $N_{\text{incre}}$ left singular vectors from the matrix of projection errors. Here, $N_{\text{incre}}$ can be chosen as a fixed value or adaptively determined based on an energy criterion in each greedy iteration. The complete POD-Greedy procedure is summarized in Algorithm \ref{alg:POD-Greedy}.

\begin{algorithm}[h]
\caption{POD-Greedy Algorithm for Linear Basis Augmentation}
\label{alg:POD-Greedy}
\begin{algorithmic}[1]
    \Require Current linear basis $\bV_r$, basis increment $N_{\text{incre}}$, high-fidelity solution $\bxi^{m+1}$
    \Ensure Updated linear reduced basis matrix $\bV_r$
    \Statex
    \State $\bE_{\text{POD}} \gets \bxi^{m+1} - \bV_r (\bV_r^\top \bxi^{m+1})$ \Comment{Compute POD projection error}
    \State Compute truncated SVD of $\bE_{\text{POD}}$, keep first $N_{\text{incre}}$ left singular vectors as $\bV_{r,\text{incre}}$
    \State $\bV_{r} \gets [\bV_r, \bV_{r,\text{incre}}]$ \Comment{Enlarge the basis}
    \State $r \gets r + N_{\text{incre}}$ \Comment{Update dimension counter}
    \State \Return $\bV_r$
\end{algorithmic}
\end{algorithm}

\begin{algorithm}[H]
\caption{Double-Greedy Algorithm for Optimal Regularization Parameter Selection}
\label{alg:regularization}
\begin{algorithmic}[1]
    \Require Current linear basis $\bV_{r}$, snapshot matrix $\bS$, matrix $\mathbf{W}^{\top}$, and projection error matrix $\bE^{\top}$ of $m+1$ representative parameters $\{\bmu^i\}_{i=1}^{m+1}$, training set $\Xi_{\mathrm{train}}$, candidate set $\Xi_{\lambda} = \{\lambda_1, \dots, \lambda_k\}$
    \Ensure Optimal regularization parameter $\lambda^\ast$, quadratic mapping matrix $\bH_r^{\lambda^\ast}$
    \For{$i = 1$ \textbf{to} $k$}
        \State Compute $(\bH_r^i, E_{\lambda_i})$ using \textbf{Algorithm~\ref{alg:tikhonov}} with inputs $(\mathbf{W}^{\top}, \bE^{\top}, \lambda_i)$ \Comment{Solve for $\bH_r^i$}
        \For{$\bmu \in \Xi_{\text{train}}$}\footnotemark
            \State Solve ROM Eq.\eqref{eq:galerkinrom} with $\bH_r^i$ to obtain reduced state $\s$ \Comment{Online computation}
            \State Compute $\Delta_{r,m}(\bmu)$ via Eq. \eqref{eq:errorestimator} \Comment{Evaluate estimator}
        \EndFor
        \State $\Delta_{r,m}^{\lambda_i} \gets \max_{\bmu \in \Xi_{\text{train}}} \Delta_{r,m}(\bmu)$ \Comment{Worst-case estimator value}
    \EndFor
    \State Select $\Xi_{\lambda}^{\text{tem}}$ containing $n_{\lambda}$ candidates with smallest $\Delta_{r,m}^{\lambda_i}$ \Comment{First-greedy}
    \For{$\lambda \in \Xi_{\lambda}^{\text{tem}}$}
        \State Collect reduced states for all $\{\bmu^i\}_{i=1}^{m+1}$ into matrix $\mathbf{S}_r^\lambda$ \Comment{Compute ROM solutions}
        \State $E_\lambda \gets \| {\bS} - (\bV_{r}\mathbf{S}_r^\lambda + \bH_r^\lambda (\mathbf{S}_r^\lambda \widetilde{\otimes} \mathbf{S}_r^\lambda)) \|_F$ \Comment{True ROM error}
    \EndFor
    \State $\lambda^\ast \gets \argmin_{\lambda \in \Xi_{\lambda}^{\text{tem}}} E_\lambda$ \Comment{Second-greedy, select optimal $\lambda$}
    \State \Return $\lambda^\ast$, $\bH_r^{\lambda^\ast}$
\end{algorithmic}
\end{algorithm}

\footnotetext{If the ROM solver diverges or produces NaN values for any $\bmu \in \Xi_{\text{train}}$ during the inner loop, the entire inner loop over $\bmu$ is immediately terminated. A large penalty value is assigned to $\Delta_{r,m}^{\lambda_i}$ in Line 7, and the algorithm proceeds to the next $\lambda_i$ in the outer loop. This early termination avoids unnecessary computations for the remaining parameters in $\Xi_{\text{train}}$ and ensures that unstable configurations are not selected in the first-greedy stage.}

\textbf{Double-greedy algorithm for regularization parameter selection:} Assume that $m+1$ representative parameters have been selected during the offline phase, the linear reduced basis matrix $\bV_{r}$ has been generated via the POD-Greedy algorithm, and the high-fidelity solutions for these selected parameters constitute the snapshot matrix $\bS$. 
To balance accuracy and stability of the ROM, an optimal regularization parameter $\lambda^\ast$ is selected from a discrete candidate set $\Xi_{\lambda}$ using a double-greedy algorithm. First, for each $\lambda$ in $\Xi_{\lambda}$, the corresponding quadratic mapping matrix is computed and the maximum value of residual-based error estimator over the whole training set is evaluated. The $n_{\lambda}$ candidate parameters that yield the smallest maximum error estimators are retained in a temporary set $\Xi_{\lambda}^{\text{tem}}$. Second, for each candidate parameter in $\Xi_{\lambda}^{\text{tem}}$, the ROM solutions corresponding to the selected representative parameters are computed, and the resulting ROM approximation errors are evaluated. The optimal regularization parameter is finally chosen as the one that minimizes these ROM errors.

This double-greedy selection criterion is particularly important when the basis number is small, because in such cases the error estimator may not be sufficiently accurate to reliably guide the choice of $\lambda$. Empirically, if only the first greedy criterion is used, the selected regularization parameter tends to be overly large, thereby diminishing the effect of the quadratic manifold and reducing the method to a linear ROM. The complete double-greedy algorithm is outlined in Algorithm \ref{alg:regularization}.

{\textbf{Full hierarchical greedy algorithm on quadratic manifold:} Integrating these two greedy algorithms, we arrive at the overall nonlinear greedy ROM on quadratic manifolds that can be summarized as follows and outlined in Algorithm \ref{alg:greedyqm}}. First, an initial parameter $\bmu^1$ is selected, and its corresponding high-fidelity solution $\bxi^1$ is stored with a temporal sampling step size $l_{\text{sam}}$. The snapshot matrix is initialized as $\bS = \bxi^1$. The linear basis matrix $\bV_{r}$ ($r =r_0$) is constructed from the first $r_0$ left singular vectors of $\bxi^1$, and the quadratic mapping matrix is generated using Algorithm \ref{alg:regularization}.
During the $m$-th greedy iteration, error estimators are computed for every parameter in the training set. The $(m+1)$-th representative parameter is chosen as the one that maximizes the error estimator defined in Eq. \eqref{eq:errorestimator}. The linear reduced basis matrix is then updated via Algorithm \ref{alg:POD-Greedy}, and the snapshot matrix is augmented as $\bS = [\bS, \bxi^{m+1}]$. Finally, the quadratic mapping matrix is recomputed using Algorithm \ref{alg:regularization}.
The greedy iteration proceeds until either the error estimator drops below a prescribed tolerance or the number of reduced basis vectors attains a predefined maximum.

\begin{algorithm}[h]
\caption{Full Hierarchical Greedy Algorithm on Quadratic Manifold}
\label{alg:greedyqm}
\begin{algorithmic}[1]
    \Require Maximum iteration $m_{\text{max}}$, initial basis size $r_0$, basis increment $N_{\text{incre}}$,
           training set $\Xi_{\mathrm{train}}$, sampling step $l_{\mathrm{sam}}$, candidate set $\Xi_{\lambda} = \{\lambda_1, \dots, \lambda_k\}$
    \Ensure Linear reduced basis matrix $\bV_r$, quadratic mapping matrix $\bH_r$
        \State $m \gets 0$, $r \gets r_0$ \Comment{Initialize counters}
    \State Select initial parameter $\bmu^{m+1}$ (e.g., randomly or at midpoint) \Comment{Initialization}
    \State $\bxi^{m+1} \gets$ high-fidelity solution for $\bmu^{m+1}$, sampled every $l_{\mathrm{sam}}$ steps
    \State $\bS \gets \bxi^{m+1}$ \Comment{Initial snapshot matrix}
    \State $\bV_r \gets$ first $r$ left singular vectors of $\bS$ \Comment{Build initial linear basis}
    \While{$m < m_{\text{max}}$}
    \State Compute matrix $\mathbf{W}^{\top}$ via Eq. \eqref{eq:WMatrix}, and projection error matrix $\bE^{\top}$ of $\{\bmu^i\}_{i=1}^{m+1}$
    \State Compute $(\lambda^\ast, \bH_r)$ using \textbf{Algorithm~\ref{alg:regularization}} with inputs ($\bV_r$, $\bS$, $\bW^\top$, $\bE^{\top}$, $\{\bmu^i\}_{i=1}^{m+1}$, $\Xi_{\lambda}$, $\Xi_{\mathrm{train}}$) \Comment{QM construction}
         \State $m \gets m + 1$
        \For{each $\bmu \in \Xi_{\text{train}}$}
            \State Solve ROM Eq.\eqref{eq:galerkinrom} with $\bV_r, \bH_r$ to get $\s$ \Comment{Online evaluation}
            \State Compute $\Delta_{r,m}(\bmu)$ via Eq.\eqref{eq:errorestimator} \Comment{Error estimator}
        \EndFor
        \State $\bmu^{m+1} \gets \argmax_{\bmu \in \Xi_{\text{train}}} \Delta_{r,m}(\bmu)$ \Comment{Greedy parameter selection}
        \State $\bxi^{m+1} \gets$ high-fidelity solution for $\bmu^{m+1}$ \Comment{Compute new snapshot}
        \State Update $\bV_r$ using \textbf{Algorithm~\ref{alg:POD-Greedy}} with inputs $(\bV_r, N_{\text{incre}}, \bxi^{m+1})$ \Comment{Enrich linear basis}
        \State $\bS \gets [\bS, \bxi^{m+1}]$ \Comment{Augment snapshot matrix}
    \EndWhile
    \State \Return $\bV_r$, $\bH_r$
\end{algorithmic}
\end{algorithm}

\section{Numerical results}
\label{sec:numerical}

In this section, we present numerical results for four parametric partial differential equations. {To measure the performance of $r$-dimensional ROMs, we calculate the average relative error, defined as}
\begin{align*}
E_r = \frac{1}{m_{\rm{test}}}\sum_{i=1}^{m_{\rm{test}}} \frac{\lVert \bu(\cdot;\bmu^i) -\widehat{\bu}(\cdot;\bmu^i)\rVert_F}{\lVert \bu(\cdot;\bmu^i)\rVert_F},
\end{align*}
where $m_{\rm{test}}$ is the number of testing parameters, and $\bu(\cdot, \bmu^i):=[\bu^1, \bu^2, \cdots, \bu^{N_T}]$ represents the full-order solution matrix for parameter $\bmu^i$.

\subsection{One-dimensional linear transport equation}

For the problem \eqref{eq:intro:transport} introduced in the introduction, we refer to the parametric setting as Case 1. Two additional parametric configurations are considered, and all three cases are described as follows.
\begin{itemize}
\item Case 1: As described in the introduction.
\item Case 2: The variation of the Gaussian pulse varies as a parameter, $\sigma \in [0.01, 0.1]$. The peak position of the initial Gaussian pulse is fixed at $x^\ast =0.5$, and the advection velocity is $c=10$. The training and testing sets are generated by uniformly partitioning the intervals $[0.01, 0.1]$ and $[0.0113, 0.0973]$ into $41$ and $5$ grids respectively.
\item Case 3: The advection velocity $c \in [1, 10]$ varies as a parameter. The variation of the Gaussian pulse is fixed at $\sigma =0.01$, and the peak position is $x^\ast =0.5$. The training and testing sets are obtained by uniformly partitioning $[1, 10]$ and $[1.013,9.973]$ into $41$ and $5$ grids respectively.
\end{itemize}

\begin{figure}[htb!]
\centering
\includegraphics[scale=0.33]{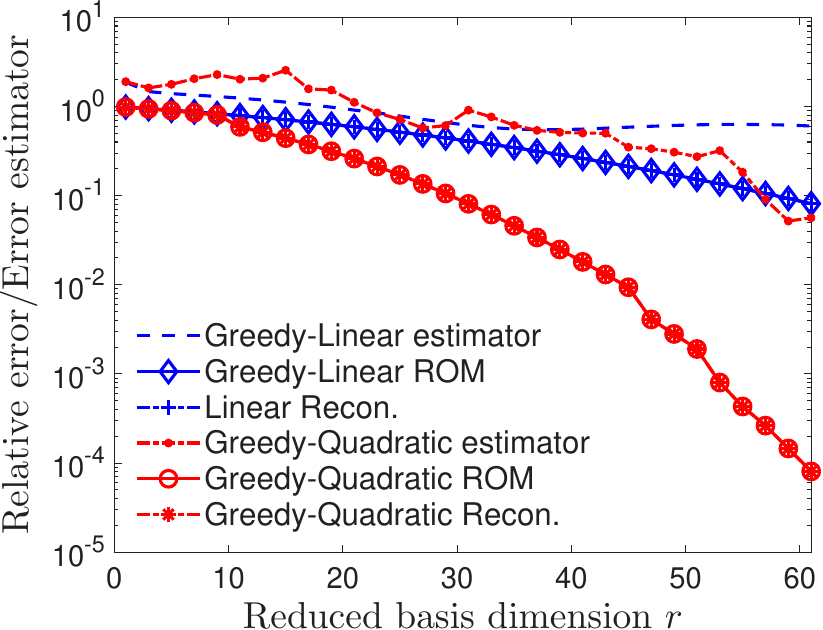}
\includegraphics[scale=0.33]{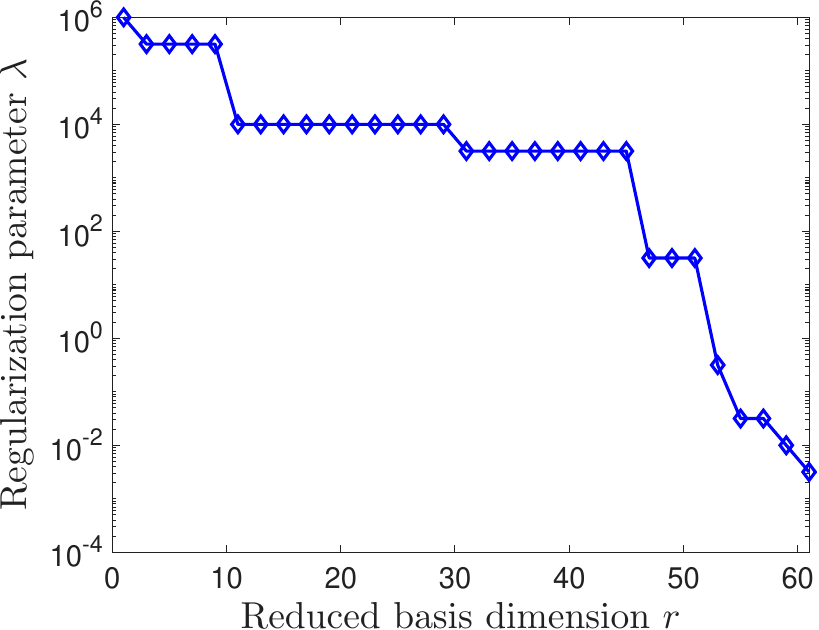}
\includegraphics[scale=0.33]{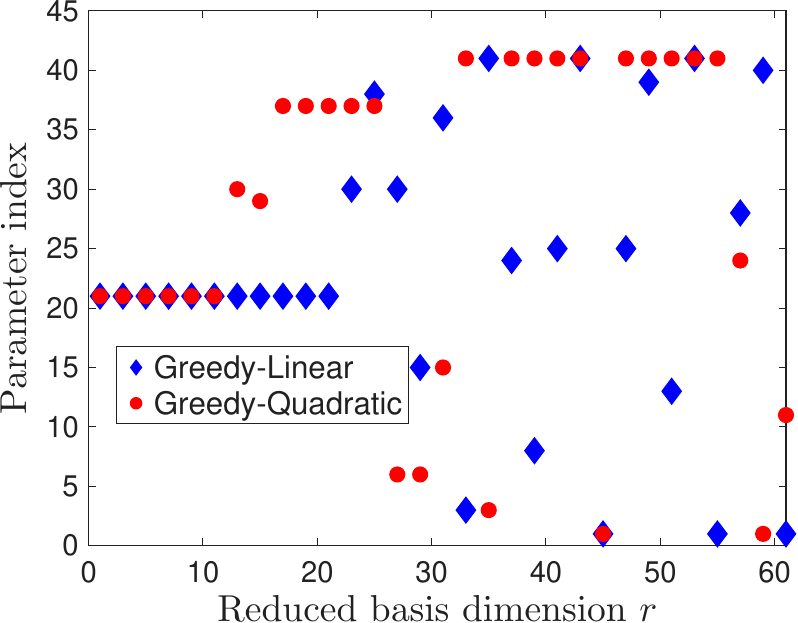}
\includegraphics[scale=0.33]{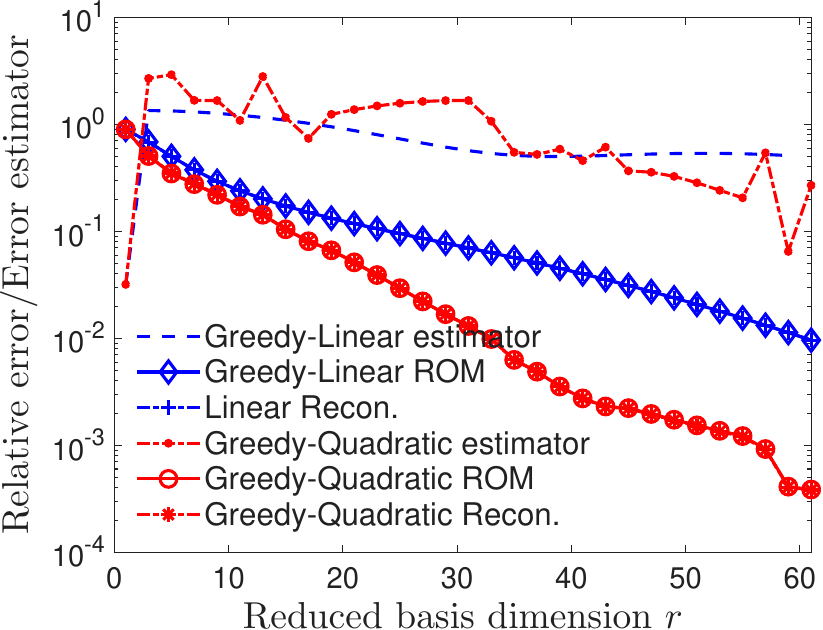}
\includegraphics[scale=0.33]{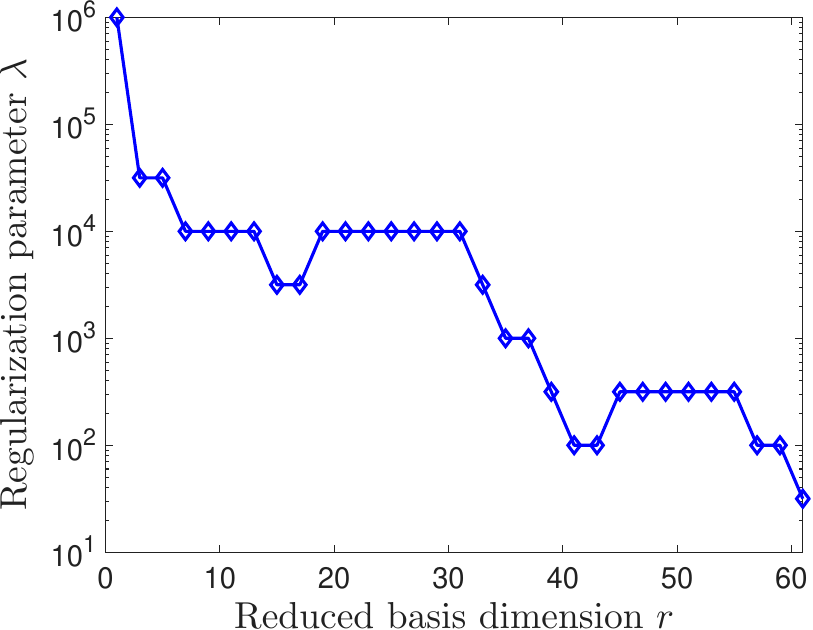}
\includegraphics[scale=0.33]{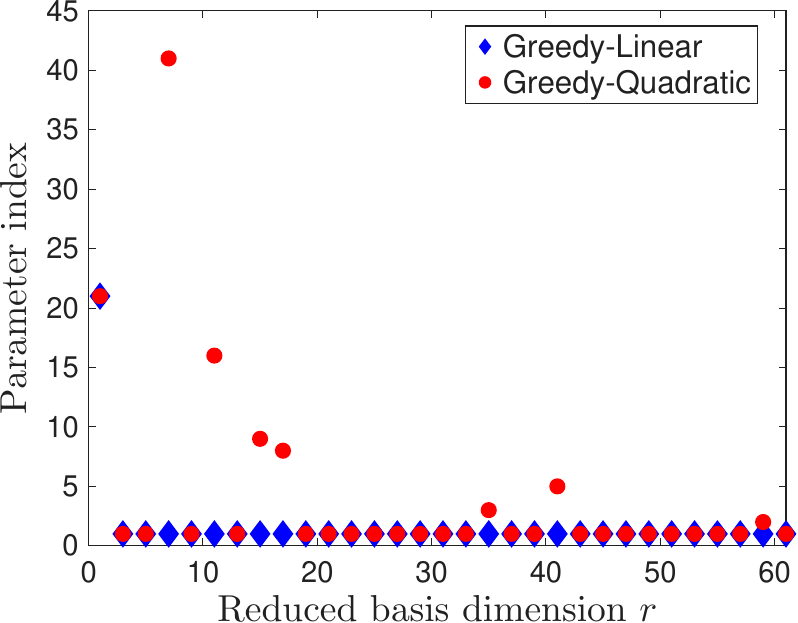}
\includegraphics[scale=0.33]{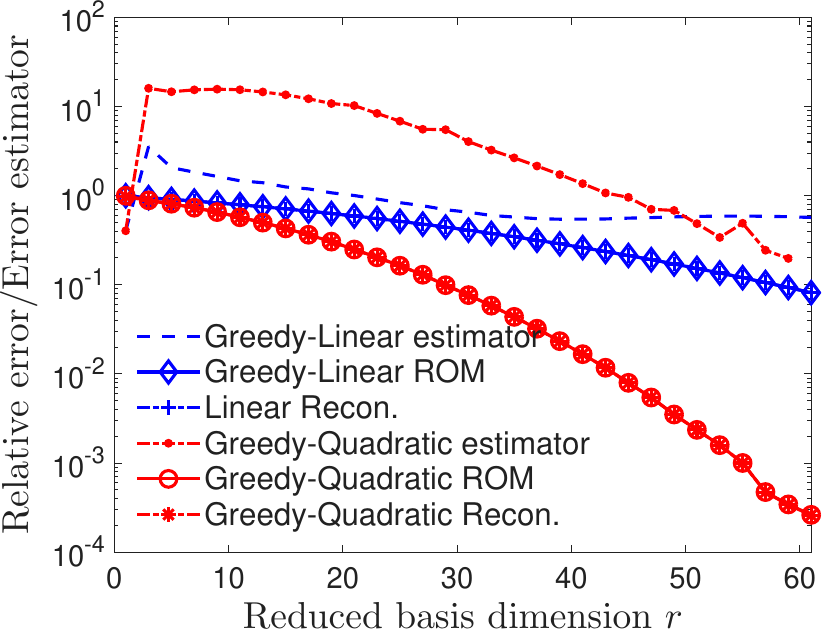}
\includegraphics[scale=0.33]{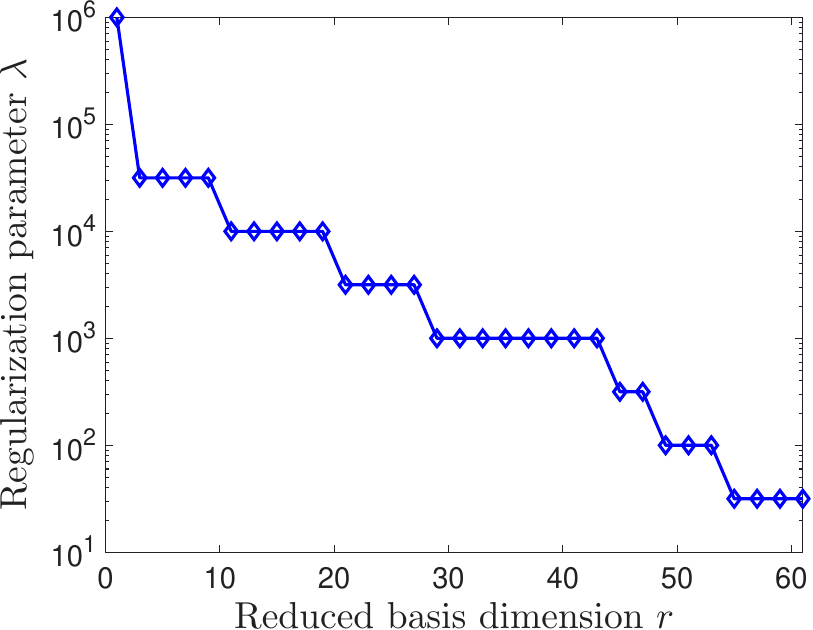}
\includegraphics[scale=0.33]{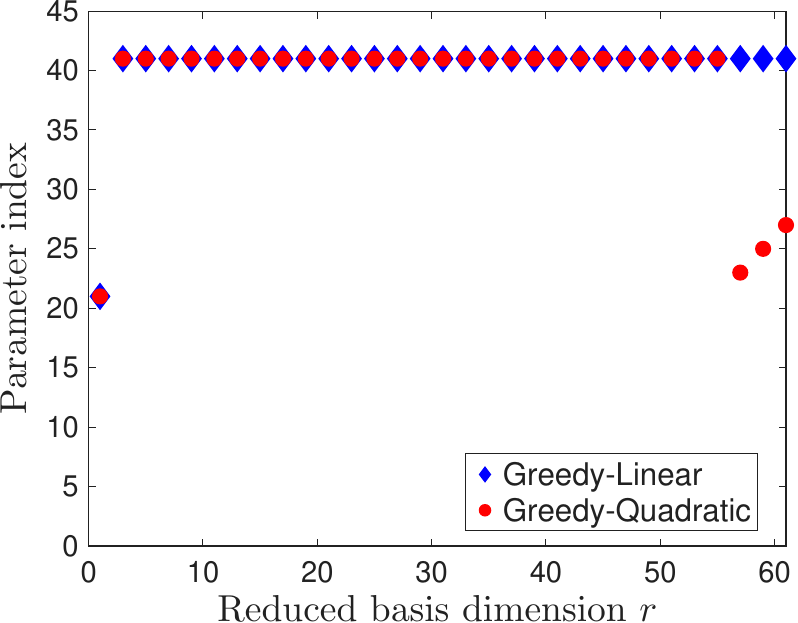}
\caption{Average relative errors, error estimators, regularization parameters, and indices of selected parameters as functions of the reduced basis dimension $r$ for $N_{\text{incre}}=2, r_0=1$, and $l_{\text{sam}} =2$. Case 1 to Case 3 (Top to Bottom row) and $n_{\lambda}=2, 3, 4$ for Case 1 to Case 3.}
\label{fig:transport:2}
\end{figure}

The candidate set for the regularization parameter in all cases is $\Xi_\lambda = 10^{-6:0.5:6}$. The resulting average relative errors and error estimators are shown in Figure \ref{fig:transport:2}. One can observe that the ROM error and reconstruction error overlap well and converge stably as the reduced dimension increases. Moreover, the residual-based \textit{a posteriori} error estimator provides an upper bound for the error curves. {For Case 1, the regularization parameters selected during the offline process are relatively large for small $r$, decrease to around $10^4$ for $10\leq r \leq 45$, and then drop sharply when $r>45$. } As the number of basis functions grows, the proposed Greedy-Quadratic method achieves progressively higher accuracy compared with the linear ROM. The representative parameters chosen by both the linear and nonlinear ROMs are distributed in a globally representative manner.
For Cases 2 and 3, the Greedy-Quadratic method likewise yields better accuracy than the linear ROM, particularly when the number of basis functions is large. The regularization parameter again decreases in a stepwise pattern. Furthermore, the distribution of selected parameters is similar for the linear and quadratic methods; specifically, parameters with smaller variance $\sigma$ (in Case 2) and larger velocity $c$ (in Case 3) are preferentially chosen during the offline process. It is worth noting that while the estimators generally decrease as $r$ increases, there are instances where the relative error is small while the estimator remains relatively large, compared to the Greedy-linear ROM. This suggests that the current estimator is not fully consistent and may require refinement in future work.

\begin{figure}[htb!]
\centering
\includegraphics[scale=0.33]{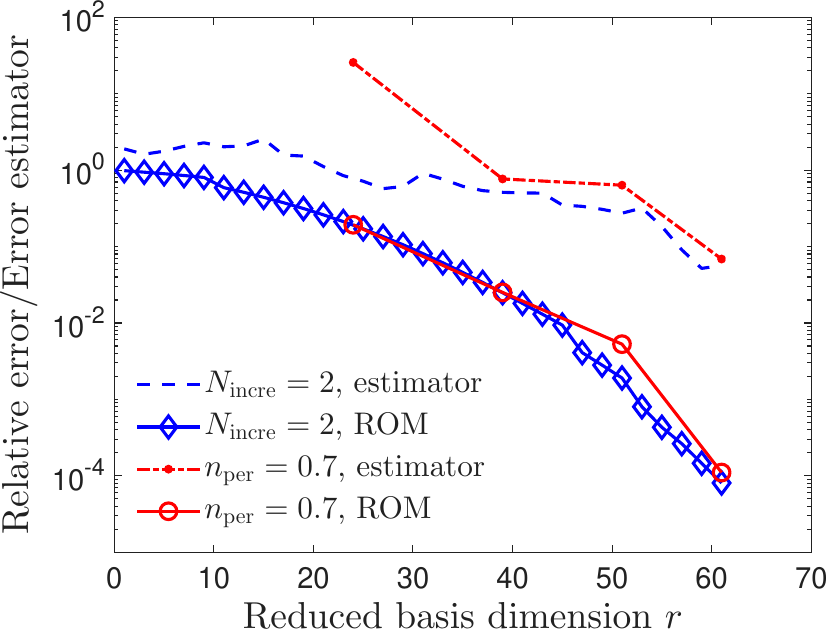}
\includegraphics[scale=0.33]{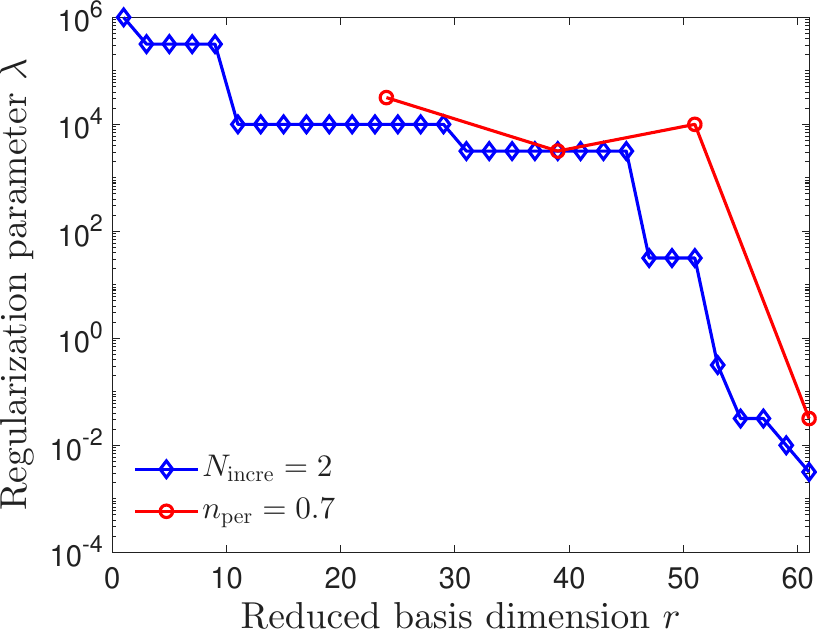}
\includegraphics[scale=0.33]{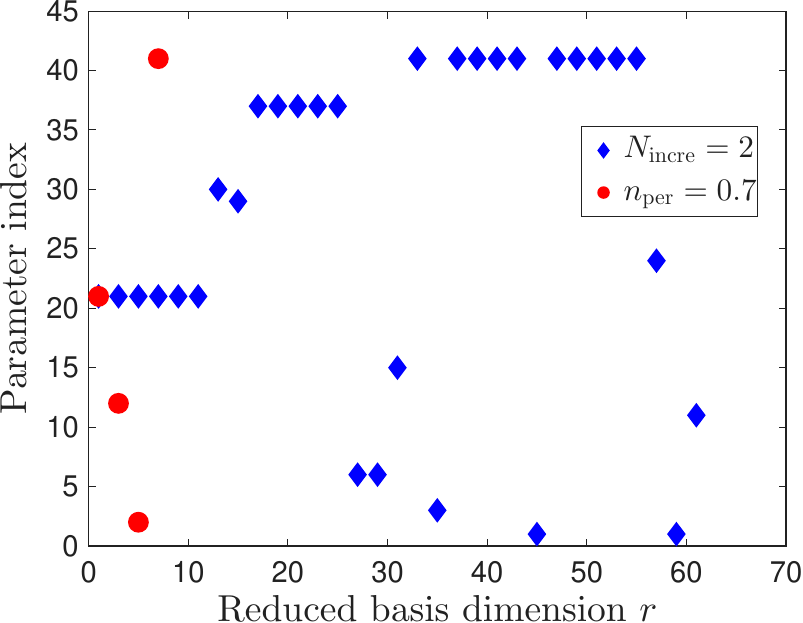}
\includegraphics[scale=0.33]{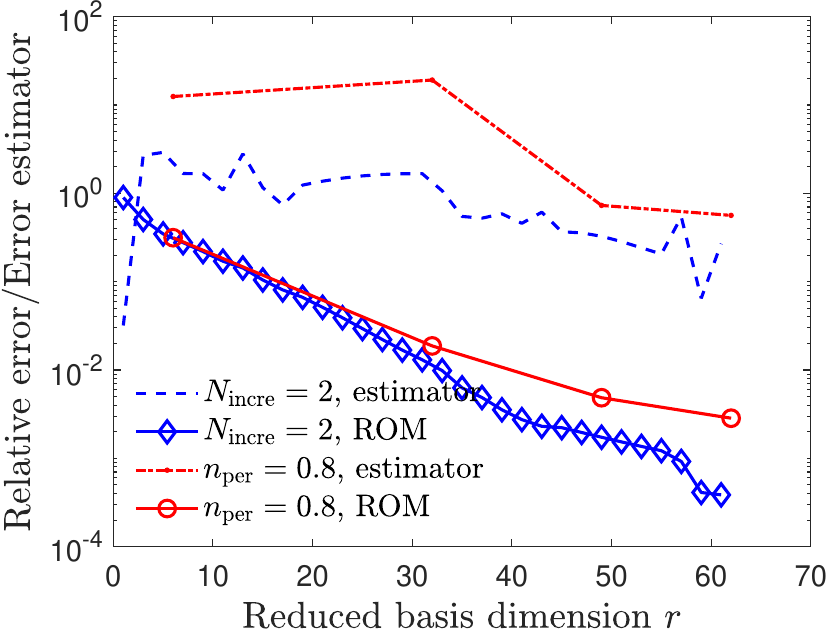}
\includegraphics[scale=0.33]{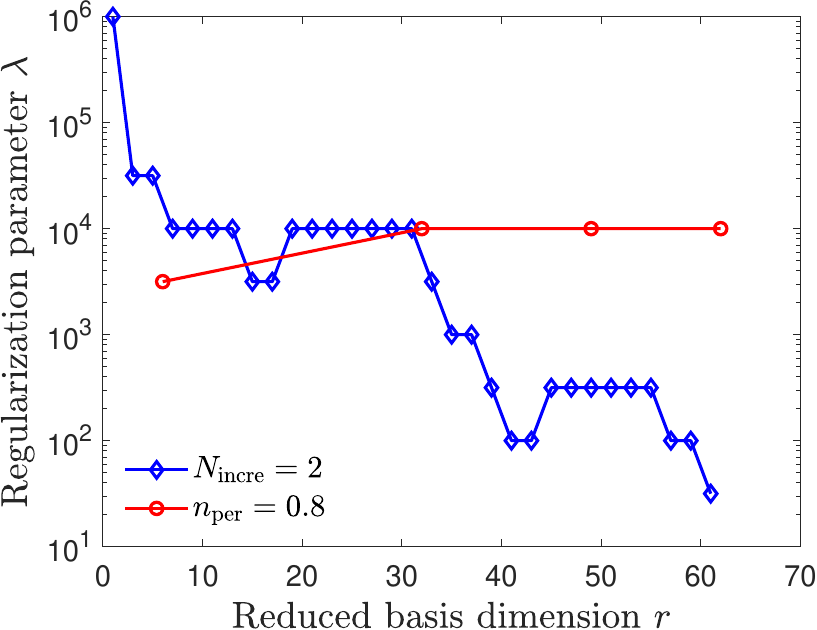}
\includegraphics[scale=0.33]{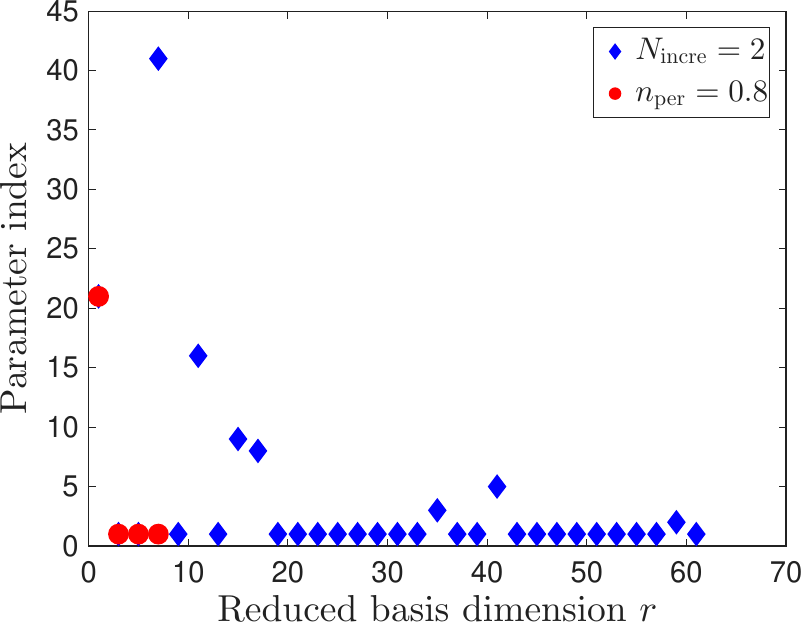}\\
\includegraphics[scale=0.33]{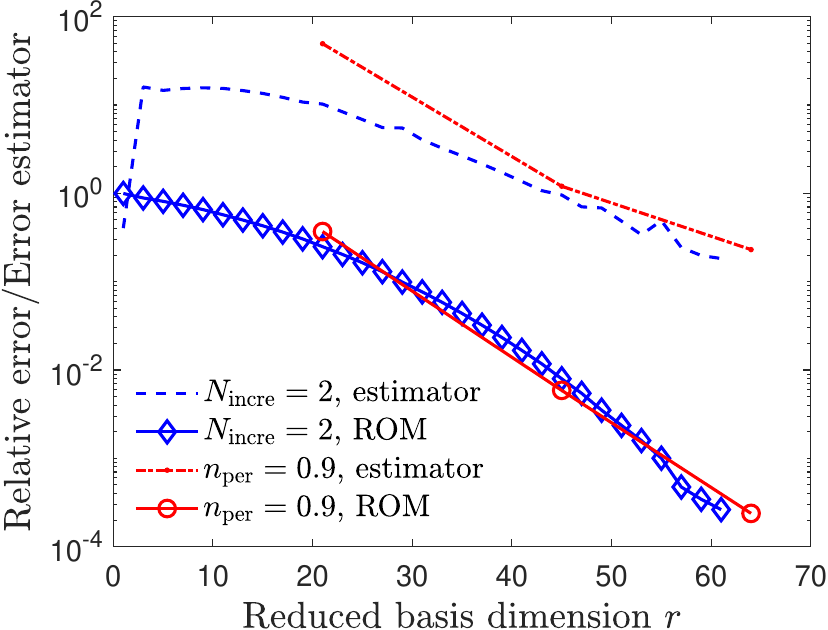}
\includegraphics[scale=0.33]{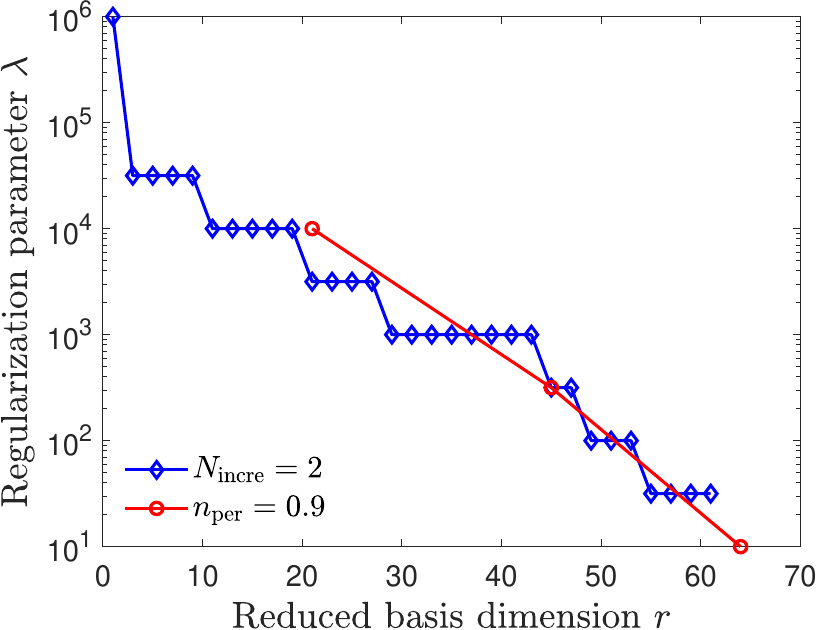}
\includegraphics[scale=0.33]{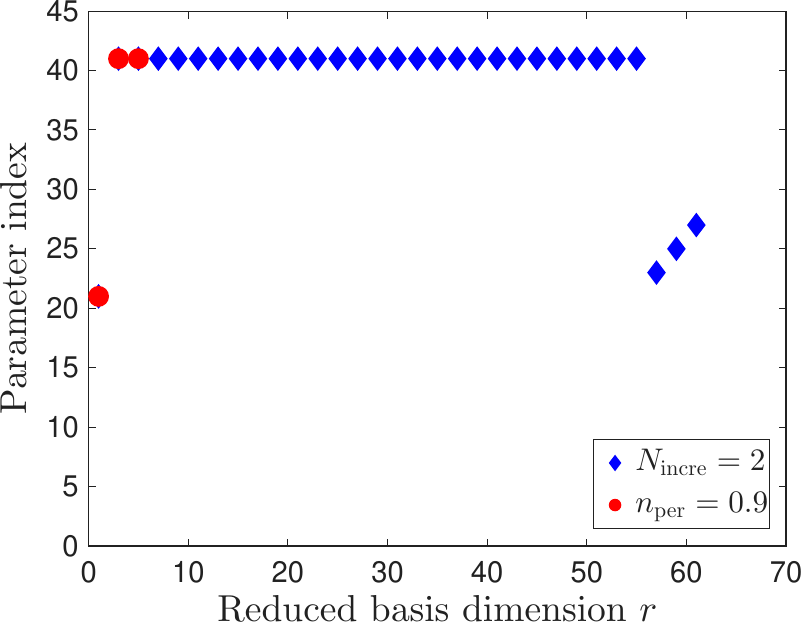}
\caption{Average relative errors, error estimators, regularization parameters, indices of selected parameters as functions of the reduced basis dimension $r$ for different $N_{\text{incre}}$. Case 1 to Case 3 (Top to Bottom row).}
\label{fig:transport:4}
\end{figure}

Figure \ref{fig:transport:2} shows that the greedy algorithm selects distinct parameters in Case~1, whereas it repeatedly selects specific parameters in Cases~2 and~3. This behavior can be explained as follows.  
In Case~1, the solution for each parameter spans a full period; hence the algorithm may choose any parameter from the set.
In Case~2, pulses with smaller variance are more difficult to approximate, making them more likely to be selected.  
In Case~3, the dynamics associated with the highest velocity are more significant for constructing the reduced basis, so they are chosen repeatedly. Moreover, because the basis increment is fixed as \(N_{\text{incre}}=2\), the reduced basis is enriched progressively, and only part of the solution space is well approximated in each greedy iteration. This explains why the same parameter may be selected in successive loops.  
If the incremental size is increased, the total number of greedy iterations decreases.  
For an intuitive illustration, a percentage-based criterion can be used to determine the incremental size: \(N_{\text{incre}}\) is chosen such that the energy captured by the POD-Greedy procedure in Algorithm \ref{alg:POD-Greedy} exceeds a given threshold $n_{\text{per}}$, where \(n_{\text{per}}=0.7, 0.8\), \text{and} $0.9$ for Case 1, Case 2, and Case 3, respectively. The corresponding results are shown in Figure \ref{fig:transport:4}. When a large percentage threshold is set, the greedy algorithm requires only two or three iterations, the total basis size reaches 64, and the accuracy remains comparable to that obtained with \(N_{\text{incre}}=2\) for Cases 1 and 3. For Case 2, when the basis number \(r > 40\), the selected regularization parameter for \(n_{\text{per}}=0.8\) is larger than that for \(N_{\text{incre}}=2\). Consequently, the average relative error achieved with \(n_{\text{per}}=0.8\) is higher. This suggests that the hyper-parameter \(n_\lambda = 3\) is unsuitable for the configuration with \(n_{\text{per}}=0.8\) and requires adjustment.
Notably, a drawback of using a large incremental size is that intermediate results for smaller basis counts are not available. 
The choice of an appropriate incremental size ultimately depends on the particular problem and the primary objectives of the model reduction. For later tests, a fixed incremental size is used to demonstrate the efficiency of our proposed greedy algorithm and double-greedy regularization technique.

\subsection{Two-dimensional linear acoustic wave equation}

The second test problem involves the two-dimensional linear acoustic wave equation with parametrized initial conditions. The governing equations are
\begin{align}
    \partial_{t} \rho(t,\boldsymbol{x}) &= -\nabla \cdot \boldsymbol{v}(t,\boldsymbol{x}), \\
    \partial_{t} \boldsymbol{v}(t,\boldsymbol{x}) &= -\nabla \rho(t,\boldsymbol{x}),
\end{align}
where $\rho(t,\boldsymbol{x})$ denotes the density field at time $t$ and position $\boldsymbol{x} = (x_1, x_2) \in [-4,4] \times [-4,4]$, and $\boldsymbol{v}(t,\boldsymbol{x})$ represents the two-dimensional velocity field. The initial density is given by a Gaussian pulse centered at $(2,2)$, parameterized by $\sigma \in [0,1]$,
\[
\rho_0(\boldsymbol{x}; \sigma) = \exp\left(-(\sigma + 6)^2 \left[(x_1 - 2)^2 + (x_2 - 2)^2\right]\right),
\]
where $\sigma$ controls the width of the pulse. The initial velocity is set to zero.
For the FOM, spatial discretization is performed using a second-order central finite difference scheme on a $500 \times 500$ grid. Time integration is carried out via a fourth-order Runge-Kutta method with a time step size of $\Delta t = T / 1200$, where $T = 6$.

\begin{figure}[htb!]
\centering
\includegraphics[scale=0.33]{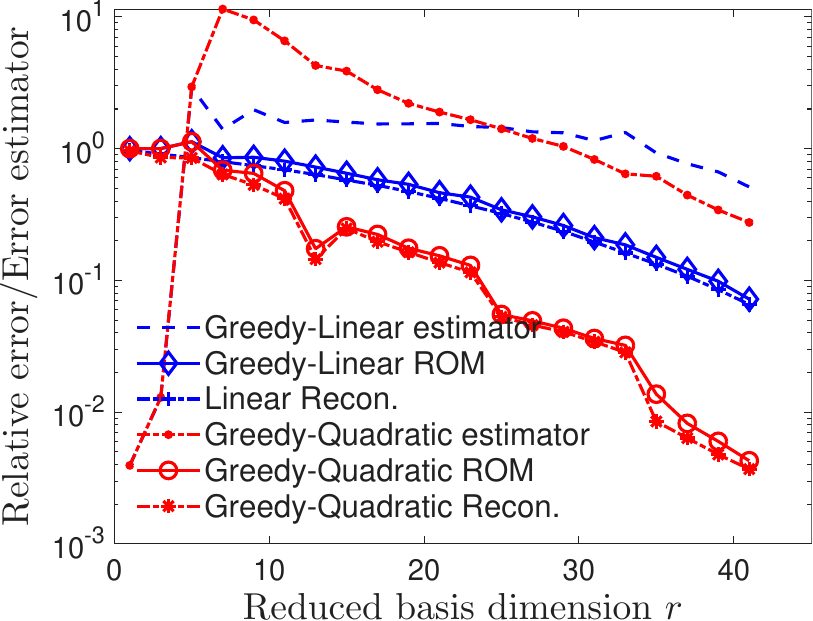}
\includegraphics[scale=0.33]{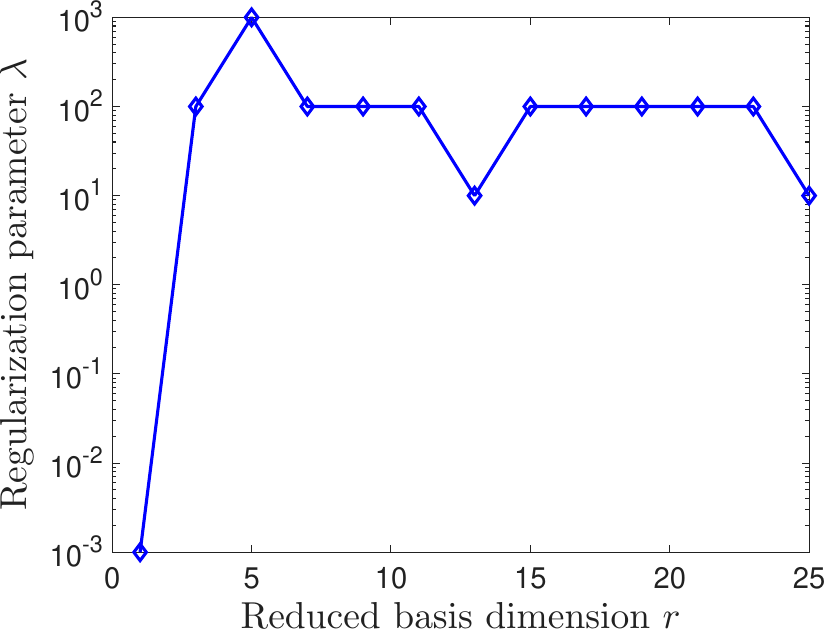}
\includegraphics[scale=0.33]{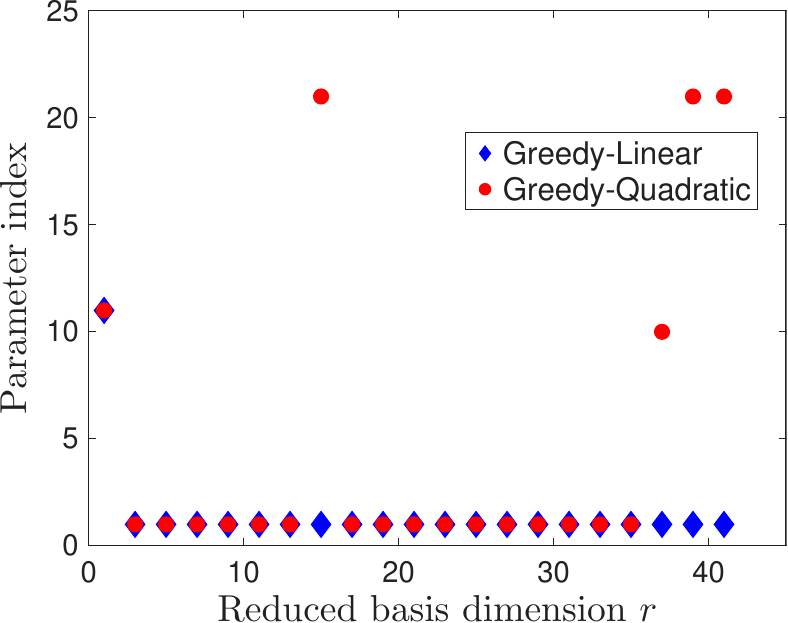}
\caption{Average relative errors, error estimators, regularization parameters, indices of selected parameters as functions of reduced basis dimension $r$ for $l_{\text{sam}}=5$, $\Xi_\lambda=10^{-3:1:3}$, and $n_{\lambda}=2$. }
\label{fig:awave:1}
\end{figure}

The training and testing sets are generated by uniformly dividing $[0,1]$ and  $[0.0053,0.953]$ into 21 and 5 grids. As shown in Figure \ref{fig:awave:1}, the average relative error of the proposed Greedy-Quadratic method is better than that of linear ROM, and converges stably when the number of reduced dimensions increases. The regularization parameters locate around $100$ and four distinct parameters $\sigma$ are selected by the greedy algorithm. Specifically, the smallest $\sigma$ is chosen at multiple greedy loops. 

\begin{figure}[htb!]
\centering
\includegraphics[scale=0.3]{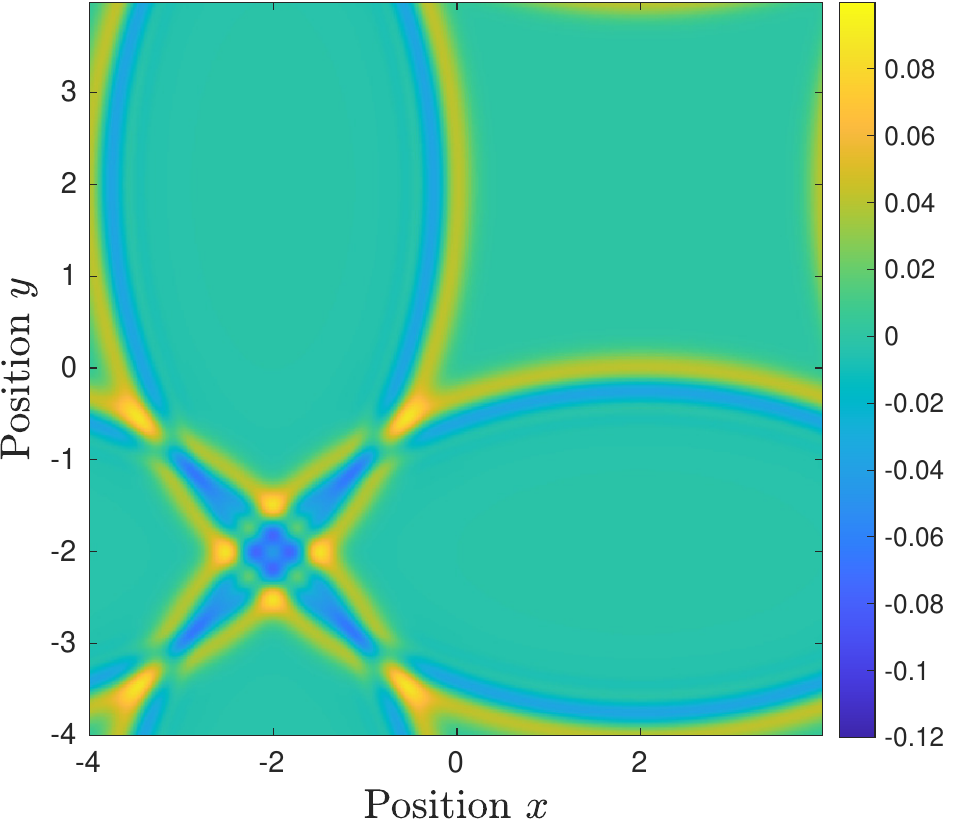}
\includegraphics[scale=0.3]{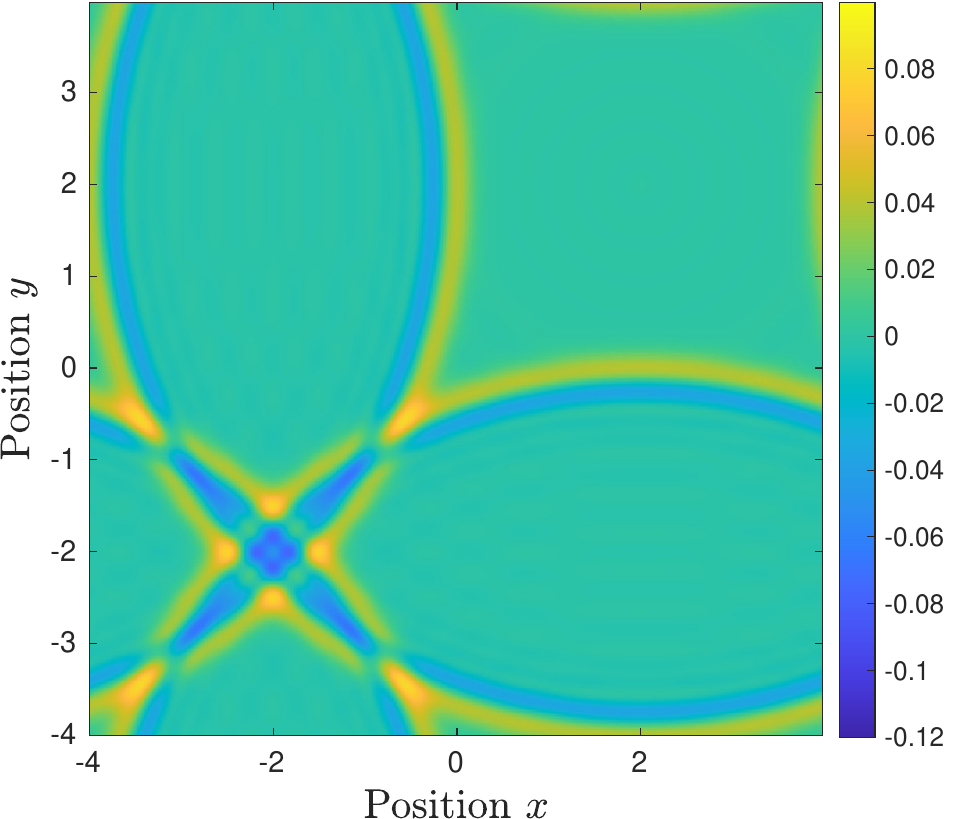}
\includegraphics[scale=0.3]{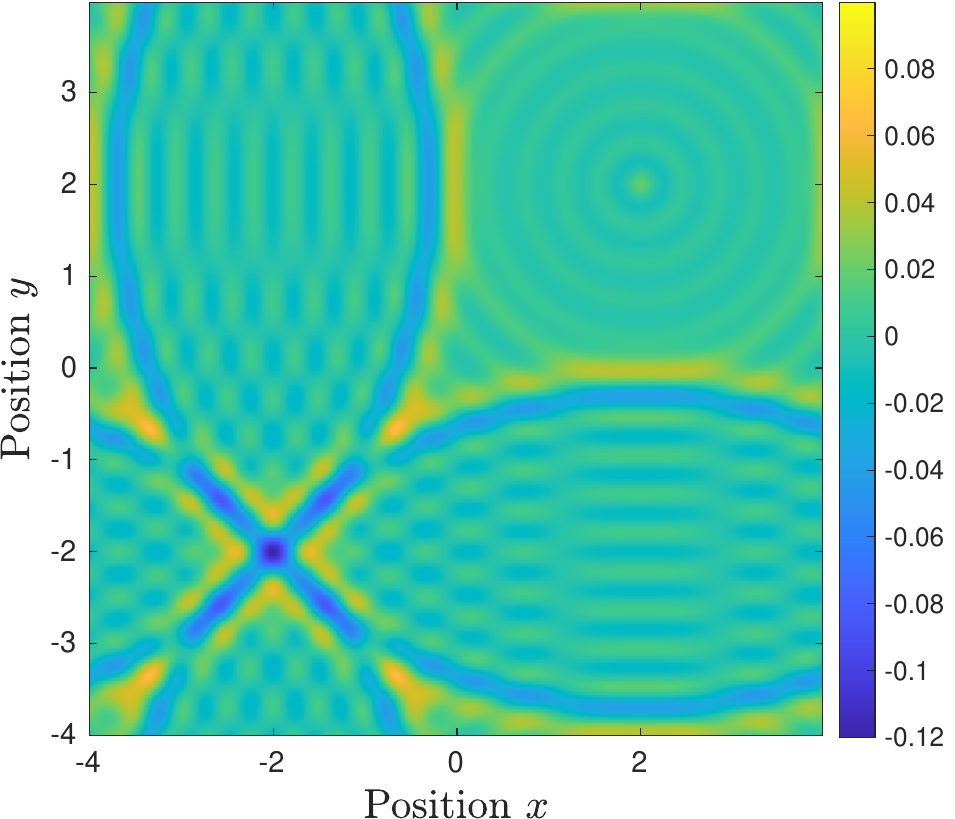}
\caption{Supremacy in accuracy of Greedy-Quadratic ROM: FDM solution, Greedy-Quadratic ROM and Greedy-Linear ROM approximations at $t=6$ with $r=29$ for $\sigma =0.47915$ (from Left to Right) for the linear acoustic wave equation.
}
\label{fig:awave:3}
\end{figure}

For an intuitive comparison, we present the approximated solutions at $t=6$ for both the Greedy-Linear and Greedy-Quadratic ROMs with $r = 29$ for $\sigma = 0.47915$ in Figure \ref{fig:awave:3}. It is evident that the solution from the Greedy-Quadratic ROM with $r = 29$ is sufficiently accurate, whereas the solution from the linear ROM still exhibits some artifacts. Regarding online acceleration, Table~\ref{tab1} lists the computational time of the high-fidelity solver as well as that of the Greedy-Quadratic ROM (with \( r = 31 \)) and the Greedy-Linear ROM (with \( r = 47 \)). At these reduced dimensions, the accuracy of the Greedy-Linear and Greedy-Quadratic ROMs is comparable. The online computational time of both ROMs is significantly lower than that of the finite difference method. The overall online computational time of the Greedy-Quadratic ROM is slightly higher, owing to the quadratic nonlinearities and the explicit time scheme.

\begin{table}[htbp]
\centering
\caption{Online time of Greedy-Quadratic ROM, Greedy-Linear ROM and FOM}
\label{tab1}
\begin{tabular}{cccc}
\toprule
Model  & Greedy-Quadratic ($r=31$)  & Greedy-Linear ($r=47$) & FOM \\
\midrule
  Time   & 0.2416& 0.0325 & 40.3881 \\
\bottomrule
\end{tabular}
\end{table}

\subsection{Two-dimensional linear advection-diffusion equation}

The third example considers a parametrized two-dimensional advection--diffusion problem \cite{haasdonk2016tutorial}, posed on the rectangular domain $\Omega = (0,2)\times(0,1)$ for $t\in[0,1]$, with the parameter vector $\bmu=(\mu_1,\mu_2)\in[0,1]^2$. The governing equation is
\begin{equation}
    \partial_t u(x,t;\bmu)
    + \nabla\!\cdot\!\big(v(x;\bmu)\,u(x,t;\bmu)\big)
    - d(\bmu)\,\Delta u(x,t;\bmu)=0,
    \label{eq:advdiff}
\end{equation}
subject to an initial condition $u(x,0;\bmu)=u_0(x;\bmu)$ and Dirichlet boundary conditions. The prescribed boundary data $g_D(x,t)$ is given by a nonnegative, time-decaying radial function centered at $x_1=1/2$, which yields a smooth inhomogeneous inflow profile.
The parametrized velocity field is
\[
v(x;\bmu)=
\begin{pmatrix}
\dfrac{\mu_1}{5}\big(1-x_2^2\big)\\[4pt]
-\dfrac12\big(4-x_1^2\big)
\end{pmatrix},
\]
and introduces a combination of horizontal shear and vertical transport. The diffusion coefficient is $d(\bmu)=0.03\,\mu_2$; thus, $\mu_2=0$ corresponds to a purely advective regime, whereas $\mu_2>0$ yields an advection--diffusion character.

For the FOM, spatial discretization is performed using a finite-volume method (FVM) on a $64\times 32$ grid. Advective terms are approximated by an explicit upwind scheme, while diffusion is treated implicitly with second-order central differences. The time interval $[0,1]$ is discretized uniformly with a time step of $\Delta t = 1/256$. Dirichlet boundary conditions are imposed through standard modifications of the discrete operators.

\begin{figure}[H]
  \centering
    \includegraphics[width=0.31\textwidth]{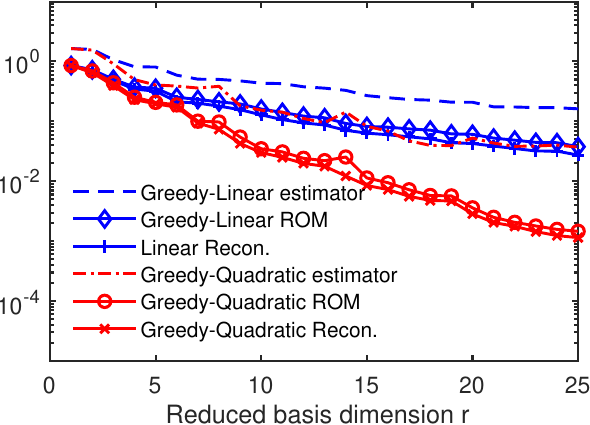}
    \includegraphics[width=0.31\textwidth]{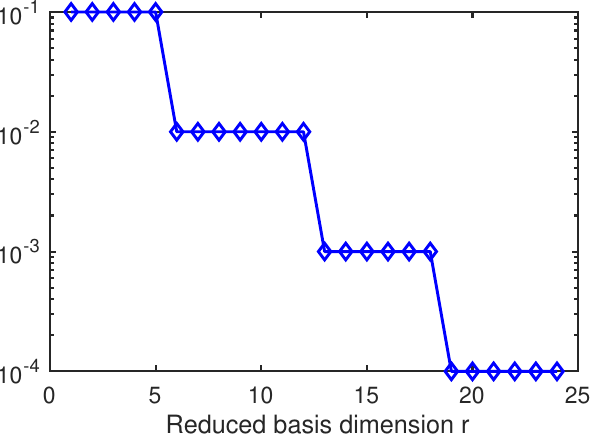} 
    \includegraphics[width=0.31\textwidth]{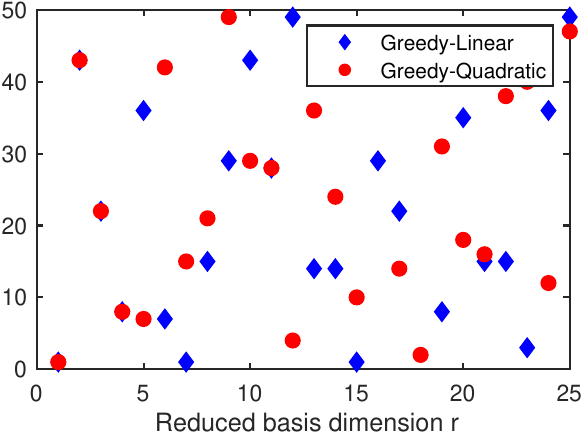} 
  \caption{Average relative errors, error estimators, regularization parameters, indices of selected parameters as functions of reduced basis dimension $r$ for $N_{\text{incre}}=1$, and $n_\lambda=1$.}
  \label{fig:alg_compare_plots}
\end{figure}

For the offline training, the training set for parameters is constructed by uniformly sampling the parameter domain $[0,1]\times[0,1]$ on a $7\times 7$ grid. The testing set is obtained by uniformly partitioning the parameter intervals $[0.19,\,0.95]$ and $[0.21,\,0.88]$ into $21$ and $5$ subintervals, respectively. The sampling time step size is $l_{\text{sam}}=1$, and the candidate regularization set is $\Xi_\lambda=10^{-5:1:0}$. As shown in Figure \ref{fig:alg_compare_plots}, the Greedy-Quadratic method achieves better accuracy as the reduced dimension increases, and the selected regularization parameter retains the stepwise decreasing feature. The selected parameters are distributed globally.

\begin{figure}[H]
\centering

% Header row
  \begin{tabular}{@{}c@{\hspace{1.2em}}c@{\hspace{1.2em}}c@{}}
    \makebox[0.30\textwidth][c]{$\bmu=(0,0)^{\mathrm{T}}$} &
    \makebox[0.30\textwidth][c]{$\bmu=(1,0)^{\mathrm{T}}$} &
    \makebox[0.30\textwidth][c]{$\bmu=(1,1)^{\mathrm{T}}$} \\
    \end{tabular}

% Row 1: FOM
\begin{tabular}{ccc}
\includegraphics[width=0.32\textwidth]{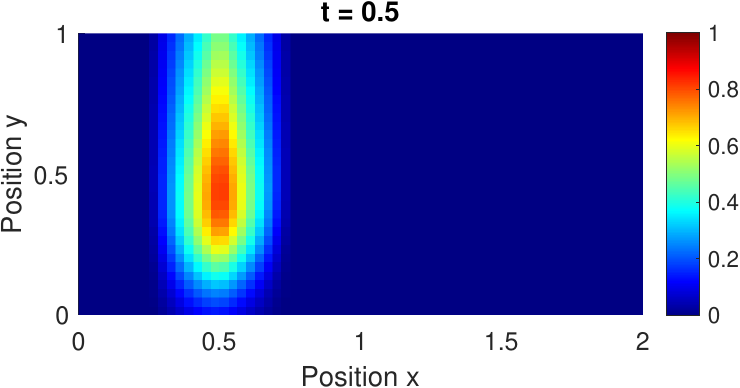} &
\includegraphics[width=0.32\textwidth]{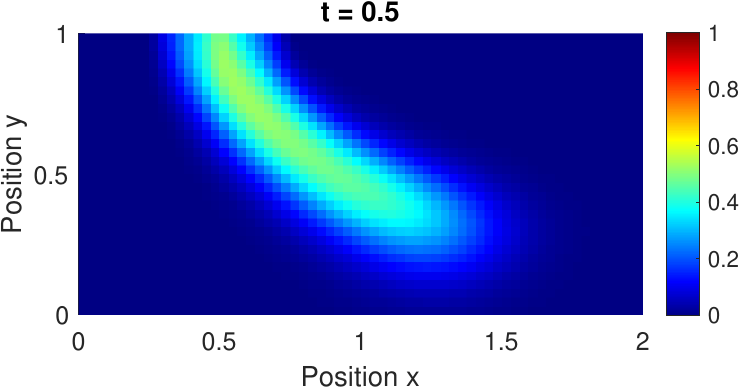} &
\includegraphics[width=0.32\textwidth]{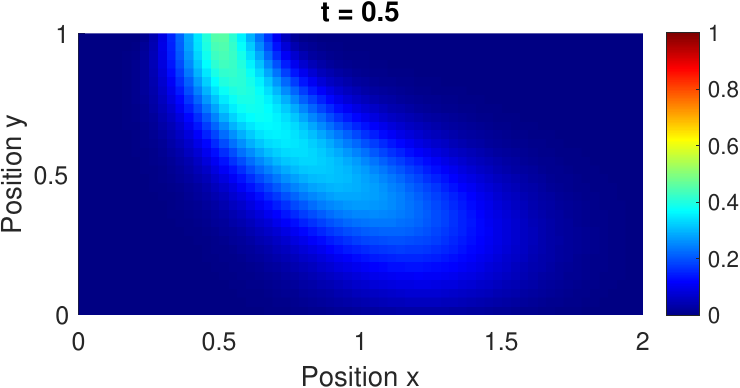} \\
\end{tabular}

\vspace{3pt}

% Row 2: Quadratic ROM
\begin{tabular}{ccc}
\includegraphics[width=0.32\textwidth]{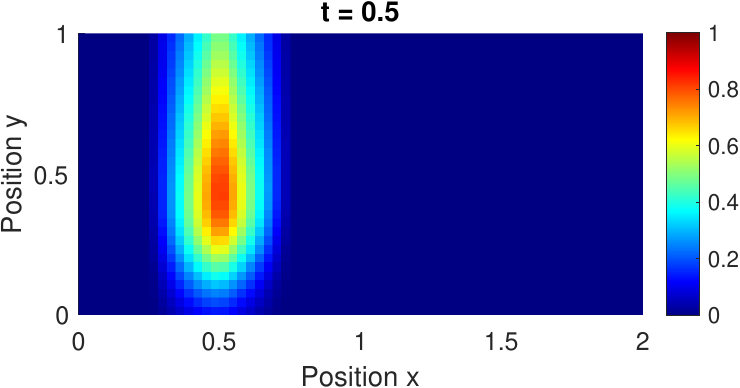} &
\includegraphics[width=0.32\textwidth]{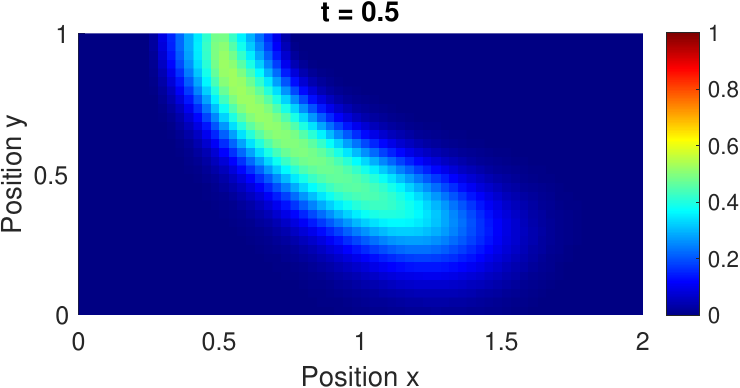} &
\includegraphics[width=0.32\textwidth]{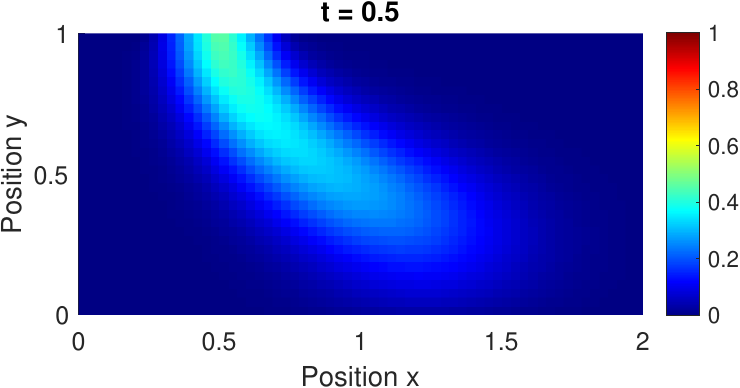} \\
\end{tabular}

\vspace{3pt}

% Row 3: Linear ROM
\begin{tabular}{ccc}
\includegraphics[width=0.32\textwidth]{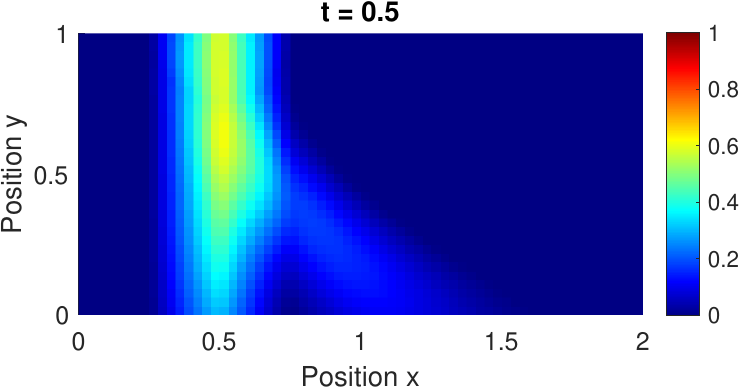} &
\includegraphics[width=0.32\textwidth]{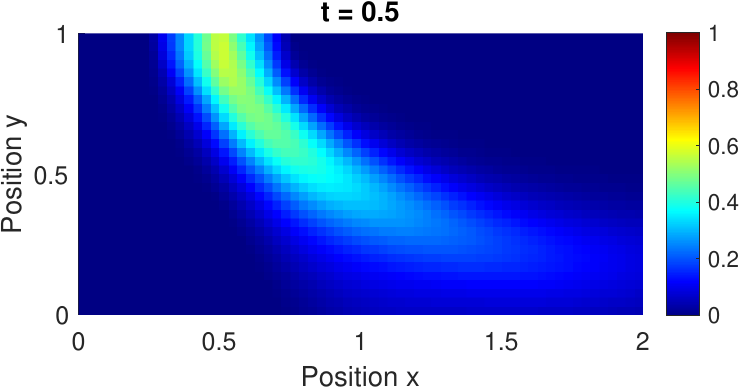} &
\includegraphics[width=0.32\textwidth]{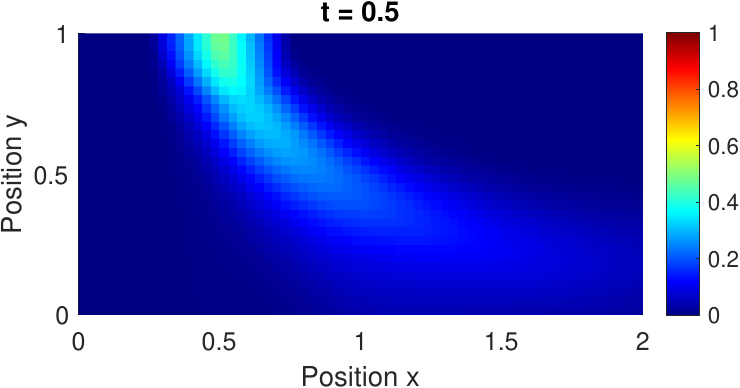} \\
\end{tabular}

    \caption{Supremacy in accuracy of Greedy-Quadratic ROM: FVM solutions, Greedy-Quadratic ROM, and Greedy-Linear ROM approximations (from top to bottom) at $t=0.5$, $r=5$ are shown for the advection diffusion equation. Columns correspond to parameter vectors $\bmu = (0,0)^{T}$, $(1,0)^{T}$, $(1,1)^{T}$ from left to right.}
    \label{fig:hf_times}

\end{figure}

Comparison results of FVM, Greedy-Quadratic ROM, and Greedy-Linear ROM approximations at $t=0.5$ for three representative parameter values,
\[
\bmu = (0,0)^{\top},\ (1,0)^{\top},\ (1,1)^{\top},
\]
 are illustrated in Figure~\ref{fig:hf_times}. The proposed Greedy-Quadratic ROM clearly outperforms the linear ROM significantly. Regarding online acceleration, the computational time of the Greedy-Quadratic ROM with $r=27$ and the Greedy-Linear ROM with $r=82$, together with the computational time of the FVM, is reported in Table~\ref{tab:online_time}. The selected reduced basis dimensions are determined by the greedy construction process and correspond to the point at which the prescribed error tolerance is reached on the training set. Both ROMs achieve substantial speedup compared to the full-order solver. However, the Greedy–Quadratic ROM requires a nonlinear solve at each time step; in this case, the smaller reduced basis dimension compensates for this additional cost, resulting in a lower overall online computational time than the Greedy–Linear ROM, even without hyper-reduction.
 
\begin{table}[htbp]
\centering
\caption{Online time of Greedy-Quadratic ROM, Greedy-Linear ROM and FOM}
\label{tab:online_time}
\begin{tabular}{cccc}
\toprule
Model 
& Greedy-Quadratic ($r=27$)  & Greedy-Linear ($r=82$)  
& FOM \\ 
\midrule
Time 
& 0.872 
& 1.455 
& 44.354\\
\bottomrule
\end{tabular}
\end{table}

\subsection{One-dimensional nonlinear viscous Burgers' equation}

The last problem is the nonlinear viscous Burgers' equation with a parametrized initial condition and periodic boundary conditions:
\[
\partial_{t}u(x,t) + u(x,t)\partial_{x}u(x,t) - \nu\partial_{xx} u(x,t) = 0,\quad x\in [-1,1], \quad t\in [0,1]
\]
where $\nu = 8\times 10^{-4}$. The parametrized initial condition is given by
\[
u_{0}(x) = 0.3\exp\left(-\sigma^{2}(x+0.5)^{2}\right) + 1,
\]
with $\sigma \in [10,15]$.
The high-fidelity solution is obtained using a finite difference scheme, discretizing in time with a fourth-order Runge-Kutta method and in space with a second order central difference under the conservative form. The simulation utilizes a time step of $\Delta t = 1/4000$ and $\mathcal{N} = 2000$ grid points.

\begin{figure}[htb!]
\centering
\includegraphics[scale=0.33]{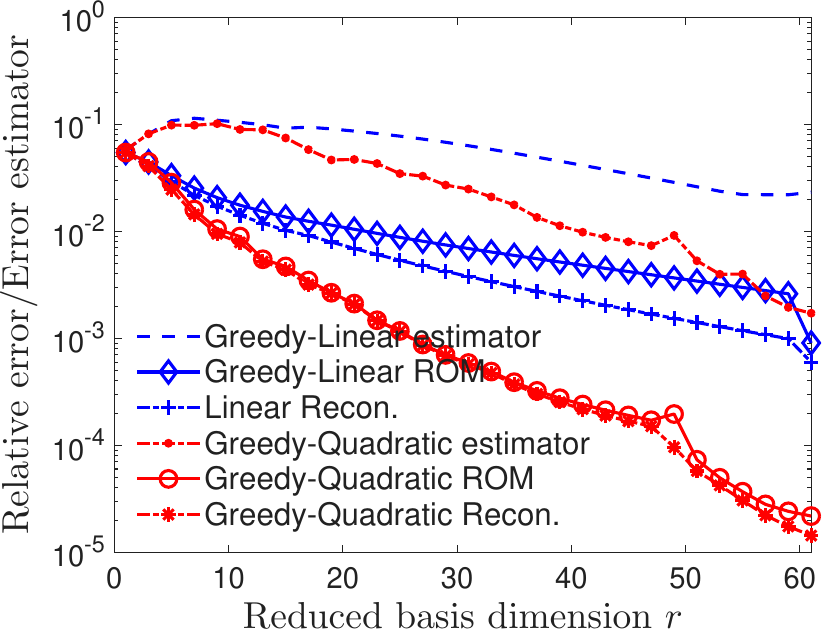}
\includegraphics[scale=0.33]{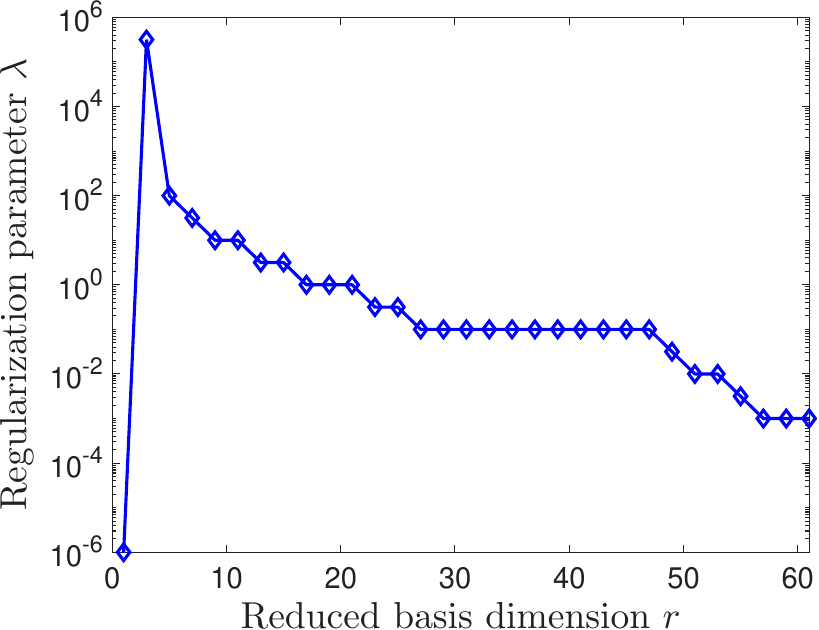}
\includegraphics[scale=0.33]{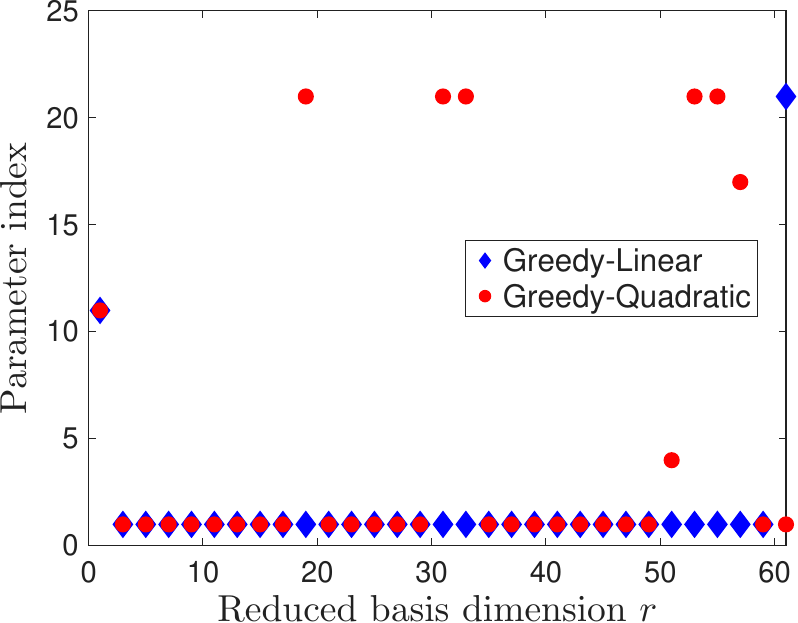}
\caption{Average relative errors, error estimators, regularization parameters, indices of selected parameters as functions of reduced basis dimension $r$ for $n_\lambda=2$.}
\label{fig:vburgers:2}
\end{figure}

\begin{figure}[htb!]
\centering
\includegraphics[scale=0.33]{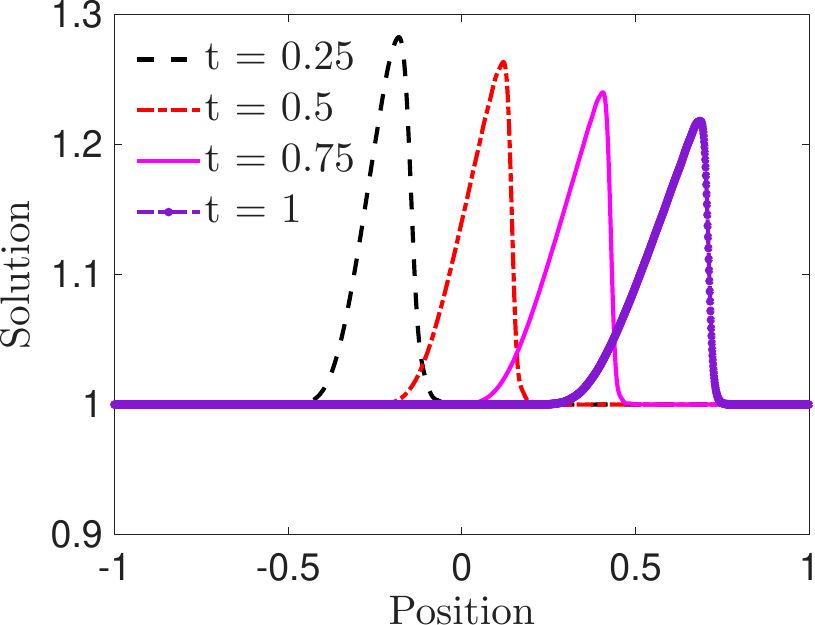}
\includegraphics[scale=0.33]{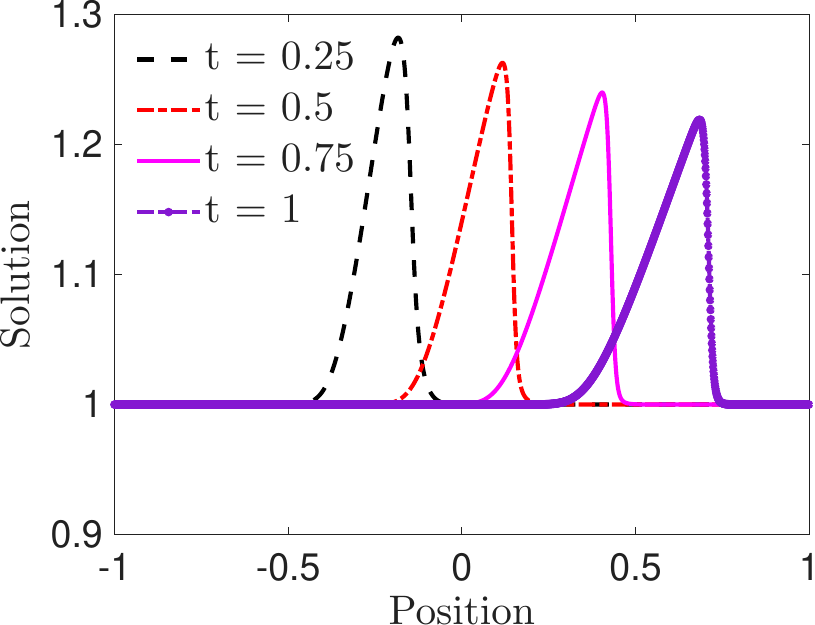}
\includegraphics[scale=0.33]{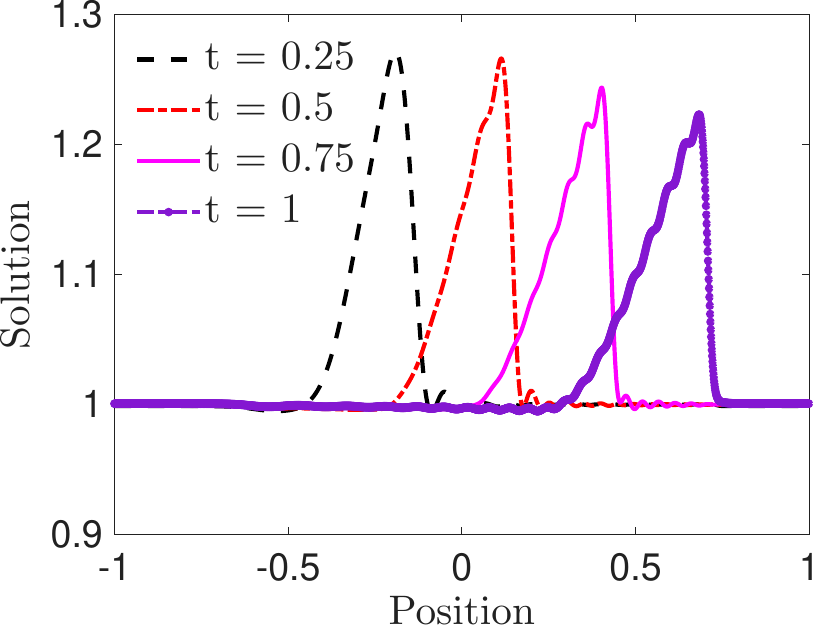}
\caption{Supremacy in accuracy of Greedy-Quadratic ROM: High-fidelity solutions, Greedy-Quadratic ROM and Greedy-Linear ROM approximations with $r=41$ for $\sigma =12.538$ at time $t=0.25, 0.5, 0.75, 1$ (from left to right) for the viscous Burgers' equation.}
\label{fig:vburgers:3}
\end{figure}

For the ROM results, the training and testing sets are generated by uniformly dividing the intervals $[10, 15]$ and $[10.123, 14.953]$ into 21 and 5 points, respectively. The time sampling step size is $l_{\text{sam}}=2$, the basis incremental is $N_{\text{incre}}=2$, and the regularization parameter set is $\Xi_\lambda= 10^{-6:0.5:6}$.
As shown in Figure \ref{fig:vburgers:2}, the error and error estimator curves for both the Greedy-Linear and Greedy-Quadratic ROMs decrease steadily. Once again, the accuracy of the Greedy-Quadratic ROM is at least one order of magnitude better than that of the Greedy-Linear ROM when the dimension of the reduced basis exceeds 30. The regularization parameter continues to exhibit a stepwise decrease. Furthermore, smaller parameter values of $\sigma$ are predominantly selected during the offline process.
For an intuitive comparison, we display the high-fidelity solutions alongside the ROM approximations at several time instants for $r=41$ in Figure \ref{fig:vburgers:3}. The solution from the Greedy-Quadratic ROM coincides with the high-fidelity solution, while the solution from the Greedy-Linear ROM exhibits noticeable oscillatory behavior.

\section{Conclusion}
\label{sec:conclusion}
This paper proposes a nonlinear ROM constructed via greedy algorithms on a quadratic manifold. An \textit{a posteriori} error estimator is formulated to guide the selection of representative parameters, thereby significantly reducing the number of queries to the high-fidelity solver. Since the error estimator provides only an upper bound on the error, it may not accurately reflect the true error when the ROM dimension is small. To address this, a double-greedy algorithm for dimension-dependent regularization ensures an appropriate balance between the accuracy and stability of the reduced system. The present work does not incorporate hyper-reduction techniques; future research may focus on the development of a more rigorously designed error estimator and the integration of hyper-reduction methods, which could further accelerate the online computation phase.

\end{document}